\documentstyle{article} 

\newtheorem{theorem}{\sc Theorem}
\newtheorem{lemma}{\sc Lemma}
\newtheorem{coro}{\sc Corollary}

\newtheorem{defin}{\sc Definition}
\newtheorem{nota}{\sc Notation}
\newtheorem{cla}{\sc Claim}
\newtheorem{ex}{\sc Example}
\newenvironment{proof}{\par \sc Proof.\rm}{\hspace*{\fill}$\Box$\vspace{1ex}}
\newenvironment{example}{\begin{ex}}{\hspace*{\fill}$\Diamond$\end{ex}}
\newenvironment{comment}{\begin{quotation}\hspace{-0.23in}\rm}{\end{quotation}}
\newenvironment{claim}{\begin{cla}}{\end{cla}}
\newenvironment{corollary}{\begin{coro}}{\end{coro}}
\newenvironment{definition}{\begin{defin}}{\end{defin}}
\newenvironment{notation}{\begin{nota}}{\end{nota}}

\newcommand{\ar}{\rightarrow}

\newcommand{\C}{C}
\newcommand{\K}{K}

\begin{document}
\title{Randomness}

\author{Paul Vit\'{a}nyi\thanks{(1944--- ) Partially
supported by the European Union
through NeuroCOLT ESPRIT Working Group Nr. 8556,
and by  NWO through NFI Project ALADDIN under Contract
number NF 62-376.
Address: CWI,
Kruislaan 413, 1098 SJ Amsterdam, The Netherlands. Email: paulv@cwi.nl}\\
CWI and Universiteit van Amsterdam}
\date{}
\maketitle

\begin{abstract}
Here we present in a single essay a combination and completion
of the several aspects
of the problem of randomness of individual objects which
of necessity occur scattered in our text~\cite{LV90}.
The reader can consult different arrangements
of parts of the material in \cite{KU87,USS90}.
\end{abstract}
\tableofcontents
\section{Introduction}
Pierre-Simon Laplace (1749---1827)
has pointed out the following reason why
intuitively a regular outcome of a random event
is unlikely.
\begin{comment}
``We arrange in our thought
all possible events in various classes;
and we regard as %
\it extraordinary %
\rm those classes
which include a very small number.
In the game of heads and tails, if head comes up a hundred
times in a row then this appears to us extraordinary,
because the almost infinite number of combinations that
can arise in a hundred throws are divided in regular
sequences, or those in which we observe
a rule that is easy to grasp, and in irregular
sequences, that are incomparably
more numerous''. [P.S. Laplace, {\em A Philosophical Essay on Probabilities},,
Dover, 1952. 
Originally published in 1819. {Translated} from 6th {French} edition. 
Pages 16-17.]
\end{comment}
If by `regularity' we mean that the complexity is
significantly less than maximal, then the number of all
regular events is small (because by simple counting 
the number of
different objects of low complexity is small).
Therefore, the event that anyone of
them occurs has small probability (in the uniform
distribution).
Yet, the classical calculus of probabilities
tells us that 100 heads are just as probable as any other sequence
of heads and tails, even though our intuition tells
us that it is less
`random' than some others. Listen to
the redoubtable Dr. Samuel Johnson (1709---1784):
\begin{comment}
``Dr. Beattie observed, as something remarkable
which had happened to him, that he chanced to see
both the No. 1 and the No. 1000, of the hackney-coaches,
the first and the last; `Why, Sir', said Johnson,
`there is an equal chance for one's seeing those
two numbers as any other two.' He was clearly right;
yet the seeing of two extremes, each of which is in
some degree more conspicuous than the rest, could
not but strike one in a stronger manner than the sight
of any other two numbers.''
[James Boswell (1740---1795), {\it Life of Johnson}, 
Oxford University Press, Oxford, UK, 1970. (Edited by R.W. Chapman,
1904 Oxford edition, as corrected by J.D. Fleeman, third edition. Originally
published in 1791.) Pages 1319-1320.]
\end{comment}
Laplace distinguishes between the
object itself and a cause of the object.
\begin{comment}
``The regular combinations occur more
rarely only because they are less numerous.
If we seek a cause wherever we perceive symmetry,
it is not that we regard the symmetrical event as less
possible than the others, but, since this event ought to
be the effect of a regular cause or that of chance,
the first of these suppositions is more probable
than the second. On a table we see letters arranged in
this order
{\tt C o n s t a n t i n o p l e},
and we judge
that this arrangement is not the result of chance,
not because it is less possible than others, for
if this word were not employed in any language we
would not suspect it came from any particular cause, but
this word being in use among us, it is incomparably
more probable that some person has thus arranged the aforesaid letters
than that this arrangement is due to chance.'' [P.S. Laplace, {\em Ibid.}]
\end{comment}
Let us try to turn Laplace's argument into a formal one.
First we introduce some notation. If $x$ is a finite binary sequence, 
then $l(x)$ denotes the {\em length} (number of
occurrences of binary digits) in $x$. For example,
$l(010)=3$.

\subsection{Occam's Razor Revisited}
Suppose we observe a binary
string $x$ of length $l(x)=n$ and want to know
whether we must attribute the occurrence
of $x$ to pure chance or to a cause. To put things in a mathematical
framework, we define {\em chance} to mean that
the literal $x$ is produced by independent tosses of a fair coin.
More subtle is the interpretation of {\em cause} as meaning
that the computer on our desk computes $x$
from a program provided by independent tosses of a fair coin.
The chance of generating $x$ literally
is about $2^{-n}$. But the chance of generating
$x$ in the form of a short program $x^*$, the cause from which
our computer computes $x$, is at least $2^{-l(x^*)}$.
In other words, if $x$ is regular, then
$l(x^*) \ll n$, and it is about
$2^{{n} - l(x^*)}$ times more likely that $x$
arose as the result of computation from some simple cause (like a short
program $x^*$) than literally by a random process.

This approach will lead to 
an objective and absolute version of the classic maxim
of William of Ockham\index{Ockham, William of} (1290? -- 1349?),
known as Occam's razor\index{Occam's razor}:
``if there are alternative explanations for a phenomenon,
then, all other things being equal,
we should select the simplest one''.
One identifies `simplicity of an object' with
`an object having a short effective description'. In other words,
\it a priori %
\rm we consider objects with short descriptions more likely
than objects with only long descriptions.
That is, objects with low complexity have
high probability while
objects with high
complexity have low probability. 

This principle is intimately
related with problems in both probability theory
and information theory. These problems as outlined below
can be interpreted as saying that the related disciplines are not
`tight' enough; they leave things unspecified which our
intuition tells us should be dealt with.

\subsection{Lacuna of Classical Probability Theory}

An adversary claims to have a true random coin and invites
us to bet on the outcome. The coin produces a hundred heads
in a row. We say that the coin cannot be fair. The adversary,
however, appeals to probabity theory which says that each sequence
of outcomes of a hundred coin flips is equally likely, $1/2^{100}$,
and one sequence had to come up. 

Probability theory gives us no basis to challenge an outcome
{\em after} it has happened. We could only exclude unfairness
in advance by putting a penalty side-bet on an outcome
of 100 heads. But what about $1010 \ldots$? What about
an initial segment of the binary expansion of $\pi$?

\begin{description}
\item[Regular sequence]
\[ \Pr(00000000000000000000000000) = \frac{1}{2^{26}} \]
\item[Regular sequence]
\[ \Pr(01000110110000010100111001) = \frac{1}{2^{26}} \]
\item[Random sequence]
\[ \Pr(10010011011000111011010000) = \frac{1}{2^{26}} \]
\end{description}

The first sequence is regular, but what is the distinction of
the second sequence and the third? The third sequence was generated
by flipping a quarter. The second sequence is very regular:
$0,1,00,01, \ldots $. The third sequence will pass
(pseudo-)randomness tests.

In fact, classical probability theory cannot express
the notion of {\em randomness of an individual sequence}.
It can only express expectations of properties of outcomes
of random processes, that is, the expectations of properties of the
total set of sequences under some distribution.

Only relatively recently, this problem has found a satisfactory
resolution by combining notions of computability 
and statistics to express the complexity of a finite object.
This complexity is the length of the shortest
binary program from which the object can be effectively
reconstructed. It may be called the {\em algorithmic
information content} of the object. This quantity turns
out to be an attribute of the object alone, and absolute (in the
technical sense of being recursively
invariant). It is the {\em Kolmogorov complexity} of the
object. 
\subsection{Lacuna of Information Theory}
In \cite{Sh48}, Claude Elwood Shannon (1916---2001) 
assigns a quantity of information or {\em entropy} to an {\em ensemble} of
possible messages. All messages in the ensemble being equally probable,
this quantity is the number of bits needed to
count all possibilities.

This expresses the fact that
each message in the ensemble can be communicated
using this number of bits.
However, it does not say
anything about the number of bits needed to convey any
individual message in the ensemble. To illustrate this,
consider the ensemble consisting of all binary strings
of length 9999999999999999.

{\sloppy
By Shannon's measure, we require
9999999999999999 bits
on the average to encode a string in such an ensemble. However, the
string consisting of 9999999999999999 1's can be encoded in about
55 bits by expressing 9999999999\ 999999 in binary and adding the
repeated pattern `1'. A requirement for this to work is
that we have agreed on an algorithm that decodes the encoded
string. We can compress the string still further when we note that
9999999999999999 equals $3^2 \times 1111111111111111$, and that
1111111111111111 consists of $2^4$ 1's.
}

Thus, we
have discovered an interesting phenomenon: the description of
some strings can be compressed considerably,
provided they exhibit enough regularity. This observation, of
course, is the basis of all systems to express very large
numbers and was exploited early on by Archimedes (287BC---212BC) in
his treatise %
\it The Sand-Reckoner%
\rm , in which he proposes
a system to name very large numbers:
\begin{comment}
``There are some, King Golon, who think that the number
of sand is infinite in multitude [...or] that no number
has been named which is great enough to exceed its multitude. [...]
But I will try to show you, by geometrical proofs,
which you will be able to follow, that, of the numbers named
by me [...] some exceed not only the mass of sand equal in
magnitude to the earth filled up in the way described,
but also that of a mass equal in magnitude to the universe.''
[Archimedes, {\em The Sand-Reckoner}, pp. 420-429 in:
{\em The World of Mathematics, Vol. 1}, J.R. Newman, Ed.,
Simon and Schuster, New York, 1956. Page 420.]
\end{comment}
However, if regularity is lacking, it becomes more cumbersome
to express large numbers. For instance, it seems easier to
compress the number `one billion,' than the number
`one billion seven hundred thirty-five million two hundred
sixty-eight thousand and three hundred ninety-four,' even though they
are of the same order of magnitude. 

The above example shows that we need too many bits to transmit
regular objects.
The converse problem, too little bits, arises as well since
Shannon's theory of information and communication
deals with the specific
technology problem of data transmission.
That is, with the information that
needs to be transmitted in order to select an object
from a previously agreed upon set of alternatives;
agreed upon by both the sender and the receiver of the message.
If we have an ensemble consisting of the {\em Odyssey} and
the sentence ``let's go drink a beer'' then we can transmit
the {\em Odyssey} using only one bit. Yet Greeks feel
that Homer's book has more information contents.
Our task is to
widen the limited set of alternatives until
it is universal. We aim at a notion of `absolute' information
of individual objects, which is the information which by
itself describes the object completely.

Formulation of these 
considerations in an objective manner leads again to the notion of 
shortest programs and Kolmogorov complexity.

\section{Randomness as Unpredictability}\label{sect.unpredictability}
What is the proper definition of a random sequence, the `lacuna in
probability theory' we have identified above? 
Let us consider how mathematicians test randomness of 
individual sequences.
To measure randomness, criteria have been developed
which certify this quality. Yet, in recognition that they do
not measure `true' randomness, we call these criteria `pseudo'
randomness tests.\index{test!pseudo randomness}
For instance, statistical survey of initial segments of the sequence
of decimal digits of $\pi$
have failed to disclose any significant deviations of randomness.
But clearly,
this sequence is so regular that it can be described by
a simple program to compute it, and this program can be
expressed in a few bits.
\begin{comment}
``Any one who considers arithmetical methods of producing random digits is,
of course, in a state of sin. For, as has been pointed out several
times, there is no such thing as a random number---there are
only methods to produce random numbers, and a strict arithmetical
procedure is of course not such a method. (It is true that
a problem we suspect of being solvable by random methods may
be solvable by some rigorously defined sequence, but this is a deeper
mathematical question than we can go into now.)'' 
[John Louis von Neumann (1903---1957), Various techniques used in connection
with random digits, {\em J. Res. Nat. Bur. Stand. Appl. Math. Series},
3(1951), pp. 36-38. Page 36. Also, {\em Collected Works, Vol. 1},
A.H. Taub, Ed., Pergamon Press, Oxford, 1963, pp. 768-770. Page 768.]
\end{comment}
This fact prompts more sophisticated definitions
of randomness. In his famous address to the International
Congress of Mathematicians in 1900,
David Hilbert (1862---1943) proposed
twenty-three mathematical problems
as a program to direct the mathematical efforts
in the twentieth century. The 6th problem asks for  "To treat
(in the same manner as geometry)
by means of axioms, those physical sciences in which
mathematics plays an important part; in the first rank
are the theory of probability ..". Thus, Hilbert views probability
theory as a physical applied theory. This raises the question
about the properties one can expect from typical outcomes
of physical random sources, which {\em a priori}
has no relation whatsoever with an axiomatic mathematical theory
of probabilities. That is, a mathematical system has no direct
relation with physical reality. To obtain
a mathematical system that is an appropriate  model of physical
phenomena one needs to identify and codify essential properties
of the phenomena under consideration by empirical observations.

Notably Richard von Mises
(1883---1953)
proposed notions that approach the very essence of true randomness
of physical phenomena.
This is related with the construction of a formal mathematical
theory of probability, to form a basis
for real applications, in the early part of
this century. While von Mises' objective was
to justify the applications to the 
real phenomena, Andrei Nikolaevitch Kolmogorov's (1903---1987)
classic 1933 treatment constructs a purely axiomatic
theory of probability on the basis of set theoretic axioms.
\begin{comment}
``This theory was so successful, that the problem of finding
the basis of real applications of the results of the mathematical
theory of probability became rather secondary to many investigators. ...
[however] the basis for the applicability of the results of the
mathematical theory of probability to real `random phenomena'
must depend in some form on the %
\it frequency\index{frequency interpretation of probability} concept of
probability%
\rm , the unavoidable nature of which has been
established by von Mises in a spirited manner.'' [A.N. Kolmogorov,
On tables of random numbers,
$Sankhy \bar a$, Series A, 25(1963), 369-376. Page 369.] 
\end{comment}
The point made is that the axioms of probability
theory are designed so that abstract
probabilities can be computed, but nothing is said about
what probability really means, or how the concept can be
applied meaningfully to the actual world.
Von Mises
analyzed this issue in detail, and suggested
that a proper definition of probability depends on
obtaining a proper definition of a random sequence.
This makes him a `frequentist'---a supporter of the
frequency theory.

The following interpretation and formulation of this theory is due
to John Edensor Littlewood (1885---1977), The dilemma of probability theory,
{\em Littlewood's Miscellany, Revised Edition}, B. Bollob\'as, Ed.,
Cambridge University Press, 1986, pp. 71-73.
The frequency theory to interpret probability
says, roughly, that if we perform an
experiment many times, then the ratio of favorable outcomes
to the total number $n$ of experiments will, %
\it with certainty%
\rm ,
tend to a limit, $p$ say, as $n \rightarrow  \infty$. This tells us something
about the %
\it meaning %
\rm of probability, namely, the measure
of the positive outcomes is $p$. But suppose we throw a coin
1000 times and wish to know what to expect. Is 1000 enough
for convergence to happen? The statement above does not say.
So we have to add something about the rate of convergence.
But we cannot assert a %
\it certainty %
\rm about a particular
number of $n$ throws, such as `the proportion of heads
will be $p \pm  \epsilon$ for large enough $n$ (with $\epsilon$
depending on $n$)'. We can at best say `the proportion
will lie between $p \pm  \epsilon$ with at least such and such
probability (depending on $\epsilon$ and $n_0$) whenever $n > n_0$'.
But now we defined probability in an obviously circular fashion.
\subsection{Von Mises' Collectives}
In 1919 von Mises\index{Mises, R. von} proposed to eliminate the problem
by simply dividing all infinite
sequences into special random sequences (called %
{\em collectives}\index{collective}),
having relative frequency\index{frequency!relative} limits,
which are the proper
subject of the calculus of probabilities and
other sequences.
He postulates the existence
of random sequences (thereby circumventing circularity)
as certified by abundant empirical evidence,
in the manner of physical laws
and derives mathematical laws of probability as a consequence.
In his view a naturally occurring sequence can be
nonrandom or unlawful in the
sense that it is not a proper collective.
\begin{comment}
Von Mises\index{Mises, R. von} views the theory of probabilities insofar
as they are numerically representable as a physical
theory of definitely observable phenomena,
repetitive or mass events, for instance, as found in
games of chance, population statistics, Brownian motion.
`Probability' is a primitive notion of the theory
comparable to those of `energy' or `mass' in other
physical theories.

Whereas energy or mass exist in
fields or material objects,
probabilities exist only in
the similarly mathematical idealization of
collectives (random sequences). All problems
of the theory of probability consist of deriving,
according to certain rules, new collectives from
given ones and calculating the distributions
of these new collectives.
The exact formulation
of the properties of the collectives is secondary
and must be based on empirical evidence. These
properties are the existence of a limiting
relative frequency and randomness.

The property
of randomness is a generalization of the
abundant experience in gambling houses,
namely, the impossibility of a successful
gambling system. Including this principle
in the foundation of probability, von Mises argues,
we proceed in the same way as the physicists did in
the case of the energy principle. Here too, the experience
of hunters of fortune is complemented by solid experience
of insurance companies and so forth.

A fundamentally different approach is to justify
{\it a posteriori} the application
of a purely mathematically constructed theory of probability, such as the
theory resulting from the Kolmogorov axioms. Suppose, we can show that
the appropriately defined
random sequences form a set of measure one, and without
exception satisfy all laws of a given
axiomatic theory of probability.
Then it appears practically
justifiable to assume that as a result of an (infinite) experiment only
random sequences appear.
\end{comment}
Von Mises'\index{Mises, R. von}
notion of infinite random sequence of 0's and 1's (collective)
essentially appeals to the idea that no gambler,
making a fixed number of wagers of `heads', at fixed odds [say $p$ versus $1-p$]
and in fixed amounts, on the flips of a coin [with bias $p$ versus $1-p$],
can have profit in the long run from
betting according to a system instead of betting at random.
Says Alonzo Church (1903--- )\index{Church, A.}:
``this definition [below] ... while clear as to general intent, is too
inexact in form to serve satisfactorily as the basis of a mathematical
theory.'' [A. Church, On the concept of a random sequence,
{\em Bull. Amer. Math. Soc.},
46(1940), pp. 130-135. Page 130.]

\begin{definition}\label{def.mises}
\rm
An infinite sequence $a_1 , a_2 , \ldots$ of 0's and 1's
is a random sequence in the special meaning
of %
\it collective\index{collective|bold} %
\index{sequence!von Mises random|see{collective}}
\rm if the following two conditions are
satisfied.
\begin{enumerate}
\item
Let $f_n$ is the number of 1's among the first
$n$ terms of the sequence. Then
\[
\lim_{{n} \rightarrow   \infty}  {f_n \over n }= p,
\]
for some $p$, $0 < p < 1$.
\item
A {\em place-selection rule} is a
partial function $\phi$,\index{place-selection rule!according to von Mises|bold}
from the
finite binary sequences to $0$ and $1$.
It takes the values 0 and 1,
for the purpose of selecting one after the other those indices $n$ for
which $\phi (a_1 a_2  \ldots a_{n-1} ) = 1$.
We require (1),
with the same limit $p$, also for every infinite subsequence
\[
a_{{n}_1} a_{{n}_2}  \ldots
\]
obtained from the sequence by some %
{\em admissible} place-selection rule.
(We have not yet formally stated which place-selection rules are
admissible.)
\end{enumerate}
\end{definition}

The existence of a relative frequency limit
is a strong assumption. Empirical evidence from long runs
of dice throws, in gambling houses, or with death statistics
in insurance mathematics, suggests that the relative frequencies
are %
\it apparently convergent%
\rm . But clearly, no empirical evidence
can be given for the existence of a definite limit
for the relative frequency. However long the test run,
in practice it will always be finite, and whatever
the apparent behavior in the observed initial segment
of the run,
it is always possible that the relative frequencies
keep oscillating forever if we continue.

The second condition ensures that no strategy using
an admissible place-selection rule can select a subsequence
which allows different odds for gambling than a subsequence
which is selected by flipping a fair coin.
For example, let a casino use a coin with probability $p=1/4$
of coming up heads and pay-off heads equal 4 times pay-off tails. This
`Law of Excluded Gambling Strategy' says that a gambler betting
in fixed amounts cannot make more profit in the long run
betting according to a system than from betting at random.
\begin{comment}
``In everyday language we call random those
phenomena where we cannot find a regularity
allowing us to predict precisely their results.
Generally speaking, there is no ground to believe that
random phenomena should possess any definite probability.
Therefore, we should distinguish between randomness
proper (as absence of any regularity) and stochastic
randomness (which is the subject of probability
theory). There emerges the problem of finding reasons for the
applicability of the mathematical theory of
probability to the real world.'' [A.N. Kolmogorov,
On logical foundations of probability theory,
Probability Theory and Mathematical Statistics,
{\em Lecture Notes in Mathematics, Vol. 1021},
K. It\^o and J.V. Prokhorov, Eds.,
Springer-Verlag, Heidelberg, 1983, pp. 1-5.
Page 1.]
\end{comment}
Intuitively, we can distinguish between sequences that are
irregular and do not satisfy the regularity implicit
in stochastic randomness, and sequences that are irregular
but do satisfy the regularities associated with
stochastic randomness. Formally, we will
distinguish the second type from the first type
by whether or not
a certain complexity measure of the initial segments
goes to a definite limit. The complexity measure
referred to is the length of the shortest description of the prefix
(in the precise sense of Kolmogorov complexity)
divided by its length. It will turn out that
almost all infinite strings are irregular of the second type
and satisfy all regularities of stochastic randomness.
\begin{comment}
``In applying probability theory we do not confine ourselves
to negating regularity, but from
the hypothesis of randomness of the
observed phenomena we draw definite
positive conclusions.'' [A.N. Kolmogorov,
Combinatorial foundations of information theory and
the calculus of probabilities, {\em Russian Mathematical Surveys},,
38:4(1983), pp. 29-40. Page 34.]
\end{comment}
Considering the sequence as fair coin tosses with $p= 1/2$,
the second condition
in Definition~\ref{def.mises} says
there is no %
\it strategy %
\rm $\phi$
(%
\it principle of excluded gambling
system\index{principle!excluded gambling system}%
\rm )
which assures a player betting at fixed odds and in
fixed amounts, on the tosses of the coin, to make infinite gain.
That is, no advantage is gained in the long run by
following some system, such as betting `head' after
each run of seven consecutive tails, or (more plausibly)
by placing the $n$th bet `head' after the appearance
of $n+7$ tails in succession. According to von Mises,
\index{Mises, R. von}
the above conditions are sufficiently familiar
and a uncontroverted empirical generalization
to serve as the basis of an applicable calculus of
probabilities.

\begin{example}
\rm
It turns out that the naive mathematical approach to a concrete
formulation,
admitting simply %
\it all %
\rm partial functions, comes to grief as follows.
Let $a = a_1 a_2  \ldots $ be any collective.
\index{collective}
Define $\phi_1$ as
$\phi_1 (a_1  \ldots a_{i-1} ) = 1$ if $a_i  = 1$, and
undefined otherwise. But then $p = 1$.
Defining $\phi_0$ by
$\phi_0 (a_1  \ldots a_{i-1} ) = b_i $, with $b_i$ the complement
of $a_i$, for all $i$, we obtain by the second
condition of Definition~\ref{def.mises} that $p = 0$.
Consequently, if we allow functions like $\phi_1$
and $\phi_0$ as strategy, then von Mises' definition
cannot be satisfied at all.
\end{example}

\subsection{Wald-Church Place Selection}
In the thirties, Abraham Wald (1902---1950) proposed to
restrict the {\em a priori} admissible $\phi$ to any
fixed countable set of functions. Then collectives do exist.
But which
countable set? In 1940,
Alonzo Church proposed to
choose a set of functions representing
`computable' strategies.
According to Church's Thesis\index{Church's Thesis},
this is precisely the set of
\it recursive functions\index{function!recursive}%
\rm .
With recursive $\phi$, not only
is the definition completely rigorous, and
random infinite sequences do exist, but moreover they are abundant
since the infinite random sequences
with $p = 1/2$ form a set of measure one. From the existence of
random sequences with probability $1/2$,
the existence of random sequences
associated with other probabilities can be derived.
Let us call sequences satisfying Definition~\ref{def.mises} with
recursive $\phi$ {\it Mises-Wald-Church random}.\label{def.MWCrandom}
\index{sequence!Mises-Wald-Church random|bold}
That is, the involved
{\em Mises-Wald-Church place-selection rules}
\index{place-selection rule!according to Mises-Wald-Church|bold}
consist of the partial recursive functions.

Appeal to a theorem by Wald\index{Wald, A.}
yields as a corollary that the set of Mises-Wald-Church random
sequences associated with any
fixed probability has the cardinality of the continuum. Moreover,
each Mises-Wald-Church random sequence qualifies as a normal number.
(A number is %
\it normal\index{number!normal|see{sequence, normal}} %
\rm 
in the sense of \'Emile F\'elix \'Edouard Justin Borel (1871---1956)
if each digit of the base,
and each block of digits of any length,
occurs with equal asymptotic frequency.) Note however, that
not every normal number is Mises-Wald-Church random.
This follows, for instance,
from Champernowne's sequence
(or number),\index{sequence!Champernowne's}
\index{Champernowne's number|see{sequence, Champernowne's}}
$$
0.1234567891011121314151617181920  \ldots
$$
due to David G. Champernowne\index{Champernowne, D.G.} (1912--- ),
which is normal in the scale of 10
and where the $i$th digit is easily calculated from $i$.
The definition of a Mises-Wald-Church random sequence implies
that its consecutive digits cannot be effectively computed.
Thus, an existence proof for Mises-Wald-Church random
sequences is necessarily
nonconstructive.

Unfortunately, the von Mises-Wald-Church definition
is not yet good enough, as was shown by Jean Ville\index{Ville, J.} in 1939.
There exist sequences that satisfy the Mises-Wald-Church
definition of randomness, with limiting relative frequency
of ones of $1/2$, but nonetheless have the
property that
$$
{f_n \over n } \geq {1 \over 2 }{\rm \ for\ all\ }n.
$$
The probability of such a sequence of outcomes in random
flips of a fair coin is zero. Intuition:
if you bet `1' all the time against such a sequence of
outcomes, then your accumulated gain is always positive!
Similarly, other properties of randomness
in probability theory such as the Law
of the Iterated Logarithm\index{Law!Iterated Logarithm}
do not follow from the Mises-Wald-Church definition.
An extensive survey on these issues (and parts of the sequel)
is given in \cite{La87}.

\section{Randomness as Incompressibility}

Above it turned out that describing `randomness' in terms
of `unpredictability' is problematic and possibly unsatisfactory.
Therefore, Kolmogorov tried another approach. The antithesis
of `randomness' is `regularity', and a finite string which is
regular can be described more shortly than giving it literally.
Consequently, a string which is `incompressible' is
`random' in this sense.
With respect to infinite binary sequences it is seductive to
call an infinite sequence `random' if all of its initial
segments are `random' in the above sense of being
`incompressible'. Let us see how this intuition can
be made formal, and whether leads to a satisfactory solution.

Intuitively,
the amount of effectively usable information in a finite string is the size
(number of binary digits or %
{\it bits})\index{bit: binary digit}
of the shortest
program that, without additional data, computes the string
and terminates.
A similar definition can be given
for infinite strings, but in this case the program
produces element after element forever. Thus,
a long sequence of 1's such as
\begin{center}
$ \overbrace{ 11111 \ldots 1}^{10,000 \mbox{ times}}  $
\end{center}
contains little information because a program of size about
$\log 10,000$ bits outputs it:
\begin{center}
\verb"for" $i := 1$ \verb"to" $10,000$\\
     \verb"print 1"
\end{center}
Likewise, the transcendental number $\pi  =  3.1415  \ldots$,
an infinite sequence of seemingly `random' decimal digits,
contains but a few bits of information. (There is a short program
that produces the consecutive digits of $\pi$ forever.)
Such a definition would appear to make
the amount of information in a string (or other object)
depend on the particular
programming language used.

Fortunately, it can be shown
that all reasonable choices of programming languages
lead to quantification of the amount of `absolute' information
in individual objects
that is invariant up to an additive constant. We call
this quantity the `Kolmogorov complexity' of the object.
If an object is regular, then it has a shorter description
than itself. We call such an object `compressible'.

More precisely, suppose we want to describe a given object by a
finite binary string. We do not care whether the object
has many descriptions; however, each description
should describe but one object.
From among all descriptions
of an object we
can take the length of the shortest description as
a measure of the object's complexity.
It is natural to call an object `simple' if it has
at least one short description, and to call it `complex'
if all of its descriptions are long.

But now we are in danger of falling in the trap
so eloquently described in the {\em Richard-Berry paradox},
\index{Paradox!Richard-Berry|bold} where
we define a natural number as
``the least natural number that cannot be described in less
than twenty words''. If this number does exist, we have just described
it in thirteen words, contradicting its definitional
statement. If such a number does not exist, then
all natural numbers can be described in less than twenty
words.
(This paradox is described in [Bertrand Russell (1872---1970)\index{Russell, B.}
and Alfred North Whitehead\index{Whitehead, A.N.},
{\it Principia Mathematica}, Oxford, 1917].
In a footnote they state that
it ``was suggested to us by Mr.\ G.G. Berry
\index{Berry, G.G.}
of the Bodleian Library''.)
We need to look very carefully at the notion
of `description'.

Assume that each description
describes at most one object. That is, there is a
specification method $D$ which associates at most one
object $x$ with a description $y$.
This means that $D$ is a function from the set of descriptions,
say $Y$, into the set of objects, say $X$.
It seems also reasonable to require that,
for each object $x$ in $X$, there is a description $y$ in $Y$
such that $D(y) = x$.
(Each object has a description.)
To make descriptions useful we like them to be finite.
This means that there are only countably many descriptions. Since there
is a description for each object, there are also only
countably many
describable objects.
How do we measure the complexity
of descriptions?

Taking our cue from the theory of computation,
we express descriptions as finite sequences of 0's and 1's.
In communication technology, if the specification method $D$ is known
to both a sender and a receiver, then a message $x$ can be transmitted
from sender to receiver by transmitting the sequence
of 0's and 1's of a description $y$ with
$D(y)=x$. The cost of this transmission is measured by the
number of occurrences of 0's and 1's in $y$, that is,
by the length of $y$. The least cost of transmission of $x$ is given
by the length of a shortest $y$ such that $D(y)=x$.
We choose this least cost of transmission
as the `descriptional' complexity of $x$ under specification
method $D$.

Obviously, this descriptional complexity of
$x$ depends crucially
on $D$.
The general principle involved is that the syntactic
framework of the description language
determines the succinctness of description.

In order to objectively compare descriptional complexities
of objects, to be able to say ``$x$ is more complex than $z$'',
the descriptional complexity of $x$
should depend on $x$ alone. This complexity can be viewed as related to
a universal description method which is {\em a priori}
assumed by all senders and receivers.
This complexity is optimal if no other description method
assigns a lower complexity to any object.

We are not really interested in optimality with respect to
all description methods.
For specifications to be useful at all it is
necessary that the mapping from $y$ to $D(y)$
can be executed in an effective manner. That is,
it can at least in principle be performed by humans or machines.
This notion has been
formalized as `partial recursive functions'. According to
generally accepted mathematical viewpoints it coincides
with the intuitive notion of effective computation.

The set of partial recursive functions
contains an optimal function which minimizes
description length of every other such function. We denote
this function by $D_0$.
Namely, for any other recursive function $D$,
for all objects $x$,
there is a description $y$ of $x$ under $D_0$ which is
shorter than any description $z$ of $x$ under $D$. (That is,
shorter up to an
additive constant which is independent of $x$.)
Complexity with respect to $D_0$ minorizes
the complexities with respect
to all partial recursive functions.

We identify the
length of the description of $x$ with respect
to a fixed specification function $D_0$ with
the `algorithmic (descriptional or Kolmogorov) complexity' of $x$.
The optimality of $D_0$ in the sense above
means that the complexity of an object $x$
is invariant (up to an additive constant
independent of $x$) under transition
from one optimal specification function to another.
Its complexity is an objective attribute
of the described object alone: it is an
intrinsic property of that object, and it does
not depend on the description formalism.
This complexity can be viewed as `absolute information content':
the amount of information which needs to be transmitted
between all senders and receivers when they communicate the
message in absence of any other {\em a priori} knowledge
which restricts the domain of the message.

Broadly
speaking, this means that all description
syntaxes which are powerful enough to express the partial
recursive functions are approximately equally succinct.
The remarkable
usefulness and inherent rightness of the theory
of Kolmogorov complexity stems from this
independence of the description method.
Thus, we have outlined the program for
a general theory of algorithmic complexity.
The four major innovations are as follows.

\begin{enumerate}
\item
In restricting
ourselves to formally effective descriptions
our definition covers every form of description
that is intuitively acceptable as being effective
according to general viewpoints in mathematics and logics.
\item
The restriction to effective descriptions
entails that there is a universal description
method that minorizes the description length or complexity
with respect to any other effective description
method. This would not be the case if we considered,
say, all noneffective description methods.
Significantly, this implies Item 3.
\item
The description length or complexity of an object
is an intrinsic attribute of the object independent
of the particular description method or formalizations
thereof.
\item
The disturbing Richard-Berry paradox above does not disappear,
but resurfaces in the form of an alternative
approach to proving Kurt G\"odel's (1906---1978) famous
\index{G\"odel, K.}
result that not every true mathematical statement
is provable in mathematics.
\end{enumerate}

\subsection{Kolmogorov Complexity}

\label{sec.complexity}
To make this treatment precise and
self-contained we briefly review notions and properties needed in the sequel.
We identify the natural numbers ${\cal N}$ and the finite binary sequences as
\[( 0, \epsilon ), (1,0), (2,1),(3,00),(4,01), \ldots ,\]
where $\epsilon$ is the empty sequence.
The {\em length} $l(x)$ 
of a natural number $x$ is the number of
bits in the corresponding binary sequence.
For instance, $l( \epsilon ) = 0$.
If $A$ is a set, then $|A|$ denotes the {\em cardinality} of $A$.
In some cases we want to encode $x$ in {\em self-delimiting}
form $x'$, in order to be able to decompose $x'y$
into $x$ and $y$. 
Short codes are obtained by iterating
the simple rule that a self-delimiting (s.d.) description of the length of $x$
followed by $x$ itself is a s.d. description of $x$. 
For example, both $x'=1^{l(x)} 0 x$ and $x''=1^{l(l(x))}0l(x)x$ are
s.d. descriptions for $x$, and 
$l(x') \leq 2l(x)+O(1)$ and $l(x'') \leq l(x) + 2l(l(x))+ O(1)$.
The string $x''$ self-delimits itself in a concatenation $x''y$
by the fact that an algorothm to retrieve $x$ works as follows.
First count the number of `1's with which $x''y$ starts out until
we find the first `0'..
This count is the length of the length of $x$, that is, 
the length of $l(x)$ which is $l(l(x))$. The binary substring of
length $l(l(x))$ following
the first $`0'$ is the binary representation of $l(x)$. The next substring
of length $l(x)$ following it is $x$ itself. So we can retrieve $x$ without
having to consider even one bit in the string after $x$. This
is why the binary code $x''$ for $x$ is called {\em self-delimiting}.

Let $\langle . \rangle : {\cal N} \times {\cal N} \ar {\cal N}$ denote a standard
computable bijective `pairing' function. Throughout this paper, we will assume
that 
\[ \langle x,y \rangle  := x''y. \]
This way, the $x''$ is has the length of $x$ but for an
additional logarithmic term.
Define $\langle x,y,z \rangle $ by $\langle x, \langle y,z \rangle  \rangle $. 

We need some notions from the theory of algorithms, see \cite{Ro67}.
Let $\phi_1 , \phi_2 , \ldots $ be a standard enumeration of the
partial recursive functions.
The (Kolmogorov) {\em  complexity} of $x \in {\cal N}$, given $y$, is defined as

\[ \C(x| y) = \min \{ l(\langle n,z \rangle ): \phi_n (\langle y,z \rangle ) = x \} .\]

This means that $\C(x|y)$ is the {\em minimal} number of bits 
in a description
from which $x$ can be effectively reconstructed, given $y$.
The unconditional complexity is defined as $\C(x)=\C(x| \epsilon )$.

An alternative definition is as follows. Let
\begin{equation}
\label{universal}
 \C_{\psi}(x|y) = \min \{ l(z): \psi (\langle y,z \rangle ) = x \}
\end{equation} 
be the conditional
complexity of $x$ given $y$ with reference to decoding function $\psi$.
Then $\C(x|y) = \C_{\psi}(x|y)$
for a universal partial recursive function $\psi$
that satisfies $\psi (\langle y,n,z \rangle )= \phi_{n} (\langle y,z \rangle )$.

We will also make use of the {\em prefix complexity} $\K(x)$, which denotes
the shortest {\em self-delimiting\/} description.
To this end, we consider so called {\em prefix\/} Turing machines,
which have only 0's and 1's on their input tape,
and thus cannot detect the end of the input.
Instead we define an input as that part of the input tape which the machine
has read when it halts. When $x \neq y$ are two such input, we clearly
have that $x$ cannot be a prefix of $y$, and hence the set of inputs forms
what is called a {\em prefix code\/}.
We define $\K(x)$ similarly as above, with reference to a universal
prefix machine that first reads $1^n0$ from the input tape and then
simulates prefix machine $n$ on the rest of the input.

A survey is \cite{LV90}.
We need the following properties. 
Throughout `$\log$' denotes the binary logarithm.
We often use $O(f(n)) = - O(f(n))$, so that 
$O(f(n))$ may denote a negative quantity.
For each $x,y\in {\cal N}$ we have 
\begin{equation}
\C(x|y) \leq l(x) + O(1).
\label{upperbound}
\end{equation}
For each $y \in {\cal N}$ there is an $x \in {\cal N}$ of length $n$
such that $\C(x|y) \geq n$. In particular, we can set $y = \epsilon$.
Such $x$'s may be called {\em random}, since they are
without regularities that can be used to compress
the description. Intuitively,
the shortest effective description of $x$
is $x$ itself. In general, for each $n$ and $y$, there are at least
$2^n -2^{n-c} +1$ distinct $x$'s of length $n$ with
\begin{equation}
\C(x|y) \geq n-c.
\label{random}
\end{equation}

In some cases we want to encode $x$ in {\em self-delimiting}
form $x'$, in order to be able to decompose $x'y$
into $x$ and $y$. 
Good upper bounds on the prefix complexity of $x$ are obtained by iterating
the simple rule that a self-delimiting (s.d.) description of the length of $x$
followed by $x$ itself is a s.d. description of $x$. 
Since $x'=1^{l(x)} 0 x$ and $x''=1^{l(l(x))}0l(x)x$ are
both s.d. descriptions for $x$, and this shows
that $\K(x) \leq 2l(x)+O(1)$ and $\K(x) \leq l(x) + 2l(l(x))+ O(1)$.

Similarly, we can encode $x$
in a self-delimiting form of its shortest program $p(x)$
($l(p (x))=\C(x)$) in $2\C(x) + 1$ bits.
Iterating this scheme, we can encode $x$ as
a selfdelimiting program of $\C(x)+2 \log \C(x) +1$ bits, 
which shows that $\K(x) \leq \C(x)+2 \log \C(x) +1$,
and so on.

If $\omega = \omega_1 \ldots \omega_k \ldots \omega_l \ldots $ 
is a finite or infinite
string then $\omega_{k:l}$ denotes $\omega_k \ldots \omega_l$,
a string of length $l-k+1$. 

\subsection{Complexity Oscillations}
Consider the question of how $C$ behaves
in terms of increasingly long initial segments of a fixed
infinite binary sequence (or string) $\omega$. For instance,
is it monotone in the sense that $C(\omega_{1:m} )  \leq C(\omega_{1:n} )$,
or $C(\omega_{1:m} | m)  \leq C(\omega_{1:n} | n)$,
for all infinite binary sequences $\omega$ and all $m  \leq n$?
We can readily show
that the answer is negative in both cases.
A similar effect arises when we try to use Kolmogorov complexity
to solve the problem of finding
a proper definition of random infinite sequences
({\em collectives})
\index{collective}
according to the task already set by von Mises\index{Mises, R. von} in 1919,
see Section~\ref{sect.unpredictability}.

Kolmogorov's intention was
to call an infinite binary sequence $\omega$ `random' if
there is a constant $c$ such that, for all $n$, the $n$-length
prefix $\omega_{1:n}$ has $C(\omega_{1:n} )   \geq   n-c$. However,
such sequences
do not exist. A simple argument shows that
that even for high-complexity sequences, with
$C( \omega_{1:n} ) \geq n -  \log n - 2 \log \log n$ for all $n$,
this results
in so-called {\it complexity oscillations}, where
\[ \frac{n- C (\omega_{1:n} )}{ \log n }  \]
oscillates between 0 and about $1$.\index{complexity oscillation|bold}

We show that the $C$ complexity
of prefixes of each infinite binary sequence drops
infinitely often unboundedly far below its own length.
Let $\omega$ be any infinite binary sequence,
and $\omega_{1:m}$ any
$m$-length prefix of $\omega$. If $\omega_{1:m}$ is the $n$th
binary string in the lexicographical order
$0, 1, 00, \ldots $, that is,
$n = \omega_{1:m}$, $m = l(n)$, then
$C( \omega_{1:n} ) = C( \omega_{m+1:n} ) +  c$, with $c$ a constant
independent of $n$ or $m$. Namely,
with $O(1)$ additional bits of information, we can trivially
reconstruct the $n$th binary string $\omega_{1:m}$ from the length $n - l(n)$
of $ \omega_{m+1:n}$. Then we find that
$C( \omega_{m+1:n} )  \leq n - l(n) +c$ for some constant $c$
independent of $n$, whence the claimed result
follows.

Our approach in this proof makes it easy to say something
about the frequency of these complexity oscillations.
\index{complexity oscillation}
Define a `wave' function $w$ by $w(1) = 2$ and $w(i) = 2^{w(i-1)}$,
then the above argument guarantees that there are
at least $k$ values $n_1 , n_2 ,  \ldots ,n_k$
less than $n = w(k)$ such that
$C( \omega_{1:{n}_i} )  \leq n_i - g(n_i ) +c$
for all $i = 1, 2,  \ldots ,k$. Obviously, this can be improved.

The upper bound on the oscillations, $C ( \omega_{1:n} ) = n+O(1)$, is reached
\index{complexity oscillation}
infinitely often for almost all high-complexity sequences.
Furthermore, the oscillations of all
high-complexity sequences stay above
$ n - \log n - 2 \log \log n$, but
dip infinitely often below $n -  \log  n$.

Due to the complexity oscillations the idea of
identifying random infinite sequences with those
such that $C(\omega_{1:n} )   \geq   n-c$, for all $n$,
is trivially infeasible. That is the bad news. In contrast,
a similar approach in Section~\ref{sect.random.finite} for finite binary
sequences turned out to work just fine. Its justification
was found in
Per Martin-L\"of's (1942--- ) important insight
that to justify any proposed definition
of randomness one has to show that the sequences,
which are random in the stated sense, satisfy the several
properties of stochasticity we know from the theory of
probability. Instead of proving
each such property separately,
one may be able to show, once and for all, that the
random sequences introduced possess, in an appropriate sense,
all possible properties of stochasticity.

In Section~\ref{sect.expect} we show how to define
randomness of infinite sequences in Martin-L\"of's
sense, which is a formal expression of the attribute
of satisfying all effectively testable
laws of randomness---and hence satisfactorily resolves the
quest for a proper definition of random sequences. Without
proof we state the relation between Martin-L\"of randomness
and high Kolmogorov complexity.
Let $\omega$ be an infinite binary sequence.

{\rm (i)} If there exists a constant $c$ such that
$C( \omega_{1:n} )   \geq   n - c$, for infinitely
many $n$, then $ \omega$ is random in the
sense of Martin-L\"of with respect to the uniform measure.

{\rm (ii)} The set of $ \omega$, for which there exists a constant $c$
and infinitely many $n$ such that $C( \omega_{1:n} )   \geq   n - c$,
has uniform measure one.

Hence, the set of random sequences
not satisfying the condition
(i)  has uniform measure zero.

The idea that random infinite sequences are those sequences
such that the complexity of the initial $n$-length segments
is at most a fixed additive constant below $n$, for all $n$, is one
of the first-rate ideas in the area of Kolmogorov
complexity. In fact, this was one of the motivations for
Kolmogorov to invent Kolmogorov complexity in the first
place. We have seen
in Section~\ref{sect.random.infinite} that this does not
work for the plain Kolmogorov complexity $C(\cdot)$,
due to the complexity oscillations. The next result is
important, and is a culmination of the theory.
For prefix complexity $K(\cdot)$ it is indeed the case
that random sequences are those sequences for which the complexity
of each initial segment is at least its length.

\begin{comment}
\rm
The history of invention of Kolmogorov complexity
is presented in detail in \cite{LV90}. Here it is important that
A.N. Kolmogorov\index{Kolmogorov, A.N.}
suggested the above relation between the complexity
of initial segments and randomness of infinite sequences
already in 1964
[see for example A.N. Kolmogorov, {\it IEEE Trans. Inform. Theory}, IT-14:5(1968), 662-664].
Similar suggestions were explicitly or implicitly made in
some form by
R. Solomonoff in 1960, \cite{So60}, and Gregory J. Chaitin (1947--- ) in
1969, \cite{Ch69}. 
This approach being incorrect using $C(\cdot)$ complexity,
P. Martin-L\"of who visited Kolmogorov in
Moscow in 1964 and 1965 already in 1965 [published
in {\it Inform. Contr.}, 9(1966), 602-619]
developed the theory as put forth in Section~\ref{sect.random.infinite}.
Nevertheless, Kolmogorov did not abandon the general outlines of his
original idea of connecting randomness of infinite sequences
with complexity, see pp. 405-406 of [Vladimir Andreevich
Uspensky (1930--- ), {\em J.\ Symb.\ Logic},
57:2(1992), 385-412].

Simultaneously and independently, Claus P. Schnorr\index{Schnorr, C.P.}
[%
\it J.\ Comput.\ Syst.\ Sci.%
\rm , 7(1973), 376-388] for the uniform distribution, and
Leonid A. Levin (1948--- )\index{Levin, L.A.} [%
\it Sov.\ Math.\ Dokl.%
\rm , 14(1973), 1413-1416] for arbitrary computable distributions,
introduced unequal versions of complexity
(respectively,
process complexity and the monotone variant\index{complexity!monotone}
of complexity ${\it Km}(x)$,
in order to
characterize randomness. They gave the first executions
of Kolmogorov's idea by showing that
an
infinite sequence $\omega$ is random iff
$|Km(\omega_{1:n})-n|=O(1)$ (Levin)
and a similar statement for process complexity (Schnorr).

In 1974 L.A. Levin in \cite{Le74}, see also P\'eter G\'acs (1947--- )
 \cite{Ga74},
and independently in 1975 G.J. Chaitin in \cite{Ch75},
introduced the prefix version of Kolmogorov complexity which
we have denoted by $K( \cdot )$.
Chaitin in his paper again
proposed
calling an infinite sequence $ \omega$
random iff $K( \omega_{1:n} )   \geq   n - O(1)$
for all $n$.
Happily, this proposal
characterizes once more precisely those
infinite binary sequences which are
random in\index{Martin-L\"of, P.}
Martin-L\"of's sense\index{sequence!random}
(without proof attributed to C.P. Schnorr\index{Schnorr, C.P.}, 1974,
in Chaitin's paper). This important result was widely circulated, but
the first {\em published proof} appears, perhaps,
as Corollary 3.2 in
[Vladimir Vyacheslavovich V'yugin (1948--- )\index{V'yugin, V.V.},
{\em Semiotika i Informatika}, 16(1981), 14-43, in Russian].
See for historical notes [A.N. Kolmogorov\index{Kolmogorov, A.N.}
and V.A. Uspensky\index{Uspensky, V.A.},
{\em SIAM J.\ Theory Probab.\ Appl.}, 32(1987), 387-412].

Another equivalent proposal in terms of
\index{sequence!Solovay random}
constructive measure theory was given
by Robert M. Solovay (1938--- ).
The fact that such different
effective formalizations of infinite random sequences
turn out to define the same mathematical
object constitutes evidence that our intuitive notion
of infinite sequences which are effectively random coincides
\index{sequence!random}
with the precise notion of Martin-L\"of random infinite sequences.
\end{comment}

\subsection{Relation with Unpredictability}
We recall von Mises'\index{Mises, R. von} classic approach to obtain
infinite random sequences $ \omega$
\index{sequence!random}
as treated in Section~\ref{sect.unpredictability}, which formed
a primary inspiration to the work reported
in this section. It is of great interest
whether one can, in his type of formulation,
capture the intuitively and mathematically
satisfying notion of infinite random sequence
in the sense of Martin-L\"of. According to von Mises
an infinite binary sequence $ \omega$ was random
(a collective\index{collective}) if:

\begin{enumerate}
\item
$ \omega$ has the property of frequency stability with limit $p$;
that is, if $f_n =  \omega_1 +  \omega_2 + \cdots + \omega_n$,
then the limit of $f_n /n$ exists and equals $p$.
\item
Any subsequence of $ \omega$ chosen according
to an admissible place-selection
rule\index{place-selection rule!according to von Mises} has frequency
stability with the same limit $p$ as in condition 1.
\end{enumerate}

One major problem was how to define
`admissible',
and one choice was to identify it with Church's notion
\index{Church, A.}
of selecting a subsequence $\zeta_1 \zeta_2  \ldots $
of $ \omega_1  \omega_2  \ldots $
by a partial recursive function $\phi$ by $\zeta_n =  \omega_m$ if
$\phi ( \omega_{1:r} ) = 0$ precisely $m-1$ times for $r   <   m-1$
and $\phi ( \omega_{1:m-1} ) = 0$. We called these $\phi$
`place-selection rules
according to Mises-Church,'\index{place-selection rule!according to Mises-Wald-C
hurch}
and the resulting sequences $\zeta$,
Mises-Wald-Church random.\index{sequence!Mises-Wald-Church random}

In Section~\ref{sect.unpredictability} we stated
that there are Mises-Wald-Church random sequences
with limiting frequency $1/2$ which
do not satisfy effectively testable properties of
randomness like the Law of the Iterated
Logarithm
or the Recurrence Property.\index{Law!Iterated Logarithm}\index{Law!Infinite Rec
urrence}
(Such properties are by definition satisfied
by sequences which are Martin-L\"of random.)
In fact, the distinction between the two is
quite large, since there are
Mises-Wald-Church collectives\index{collective} $ \omega$ with limiting
frequency $1/2$ such that $C( \omega_{1:n} ) = O(f(n)  \log n)$
for every unbounded, nondecreasing, total recursive function $f$.
\label{comment.miseswaldchurchrandom}
Such collectives are very nonrandom
sequences from the viewpoint of
Martin-L\"of randomness where $C( \omega_{1:n} )$
is required to be asymptotic to $n$.
See Robert P. Daley (1944--- ),
\it Math.\ Syst.\ Theory%
\rm ,
9(1975), 83-94. It is interesting to point out
that although a Mises-Wald-Church random sequence
may have very low Kolmogorov complexity, it in fact has
very high time-bounded Kolmogorov complexity.
If we consider also sequences with limiting frequencies
different from $1/2$, then it becomes even easier
to find sequences which are random according to Mises-Wald-Church,
but not according to Martin-L\"of.\index{Martin-L\"of, P.}\index{sequence!random
}
Namely, %
\it any %
\rm sequence $ \omega$ with
limiting relative frequency $p$ has complexity
$C( \omega_{1:n} )  \leq H n + o(n)$ for
$H = - (p \log p + (1-p) \log (1-p))$ ($H$ is Shannon's entropy\index{entropy}).
This means that for each $\epsilon    >   0$ there
are Mises-Wald-Church random sequences $ \omega$ with
$C( \omega_{1:n} )   <   \epsilon n$ for all but finitely many $n$.

On the other hand, clearly all Martin-L\"of random sequences
are also Mises-Wald-Church random (each
admissible selection rule is an effective sequential test).

\subsection{Kolmogorov-Loveland Place Selection}
This suggests that we have to liberate our notion of
admissible selection rule somewhat in order to capture the proper
notion of an infinite random sequence using von Mises'
approach. A proposal in this direction was given
by A.N. Kolmogorov [$Sankhy \bar a$, Ser. A, 25(1963),
369-376] and Donald William Loveland (1934--- )
\index{Kolmogorov, A.N.}\index{Loveland, D.W.}
[%
\it Z. Math.\ Logik Grundl.\ Math.\ %
\rm 12%
\rm (1966), 279-294].

A Kolmogorov-Loveland admissible selection
function\index{place-selection rule!according to Kolmogorov-Loveland} to select
an
infinite subsequence
$\zeta_1 \zeta_2  \ldots $ from $ \omega =  \omega_1  \omega_2  \ldots $
is any one-one recursive
function $\phi :   \{  0, 1  \}^* \rightarrow {\cal N} \times   \{  0, 1  \}  $
from binary strings to (index, bit) pairs.
The index gives the next position in $\omega$ to be scanned,
and the bit indicates if the scanned bit of $\omega$ must
be included in $\zeta_1 \zeta_2  \ldots $.
More precisely,
\[ \phi ( z) = \left\{ \begin{array}{ll}
(i_1 , a_1 ) & \mbox{if $z= \epsilon$,} \\
(i_m , a_m ) & \mbox{if $z=z_1 z_2  \ldots  z_{m-1}$ with
$z_j = \omega_{{i}_j},
1  \leq j   <   m$.}
\end{array}
\right. \]
The sequence $\zeta_1 \zeta_2 \ldots$ is called a
{\em Kolmogorov-Loveland random sequence}.
Because $\phi$ is one-one, it must hold
that $i_m  \neq  i_1 , i_2 , \ldots ,i_{m-1}$.
The described
process yields a sequence of $\phi$-values
$(z_1 , a_1 ),  (z_2 ,   a_2 ), \ldots $.
The selected subsequence $\zeta_1 \zeta_2  \ldots $ consists
of the ordered sequence of those $z_i$'s of which the
associated $a_i$'s equal ones.

As compared to the Mises-Wald-Church approach,
the liberation is contained in the fact that
the order of succession of the terms in the
subsequence chosen is not necessarily the same as
that of the original sequence. In comparison,
it is not obvious (and it seems to be unknown)
whether a subsequence $\zeta_1 \zeta_2  \ldots $
selected from a Kolmogorov-Loveland random sequence
$\omega_1 \omega_2  \ldots $ by a Kolmogorov-Loveland
place-selection rule is itself a Kolmogorov-Loveland
random sequence. Note that the analogous property
necessarily holds for Mises-Wald-Church random sequences.
\index{sequence!Mises-Wald-Church random}

The set of %
Kolmogorov-Loveland random %
\rm sequences
is contained in the set of Mises-Wald-Church random
sequences and contains the set of Martin-L\"of random
sequences. If $\omega_1 \omega_2  \ldots $ is
Kolmogorov-Loveland random then clearly $\zeta_1 \zeta_2  \ldots $,
defined by $\zeta_i = \omega_{{\sigma} (i)}$ with $\sigma$
being a recursive permutation, is also Kolmogorov-Loveland random.
The Mises-Wald-Church notion of randomness does not
have this important property of randomness of staying
invariant under recursive permutation. Loveland
gave the required counterexample in the cited reference.
Hence, the containment of the set of Kolmogorov-Loveland random
sequences in the set of Mises-Wald-Church random sequences
is proper.

This leaves the question of whether
the containment of the set of Martin-L\"of random sequences
\index{sequence!random}
in the set of Kolmogorov-Loveland random sequences is proper.
Kolmogorov has suggested in
\index{Kolmogorov, A.N.}
[%
\it Problems Inform.\ Transmission%
\rm , 5%
\rm (1969),
3-4] without proof that there is a
Kolmogorov-Loveland random sequences $ \omega$ such that $C( \omega_{1:n} ) = O(
 \log n)$.
But Andrei Albertovich Muchnik (1958--- )\index{Muchnik, An.A.}
(not to be confused with A.A. Muchnik\index{Muchnik, A.A.})
showed that this is false since
no $ \omega$ with $C( \omega_{1:n} )  \leq cn + O(1)$ for a constant
$c   <   1$ can be Kolmogorov-Loveland random. Nonetheless,
containment is proper since Alexander Khanevich Shen' (1958--- )
\index{Shen', A.Kh.}
[%
\it Soviet Math.\ Dokl.%
\rm , 38:2(1989), 316-319]
has shown there exists a Kolmogorov-Loveland random sequence
\index{sequence!Kolmogorov-Loveland random}
which is not random in Martin-L\"of's sense.
Therefore, the problem of giving a satisfactory
definition of infinite Martin-L\"of random sequences in
the form proposed by von Mises has not yet
been solved. See also [A.N. Kolmogorov\index{Kolmogorov, A.N.}
and V.A. Uspensky,\index{Uspensky, V.A.} {\it Theory Probab.\ Appl.},
32(1987), 389-412; V.A. Uspensky, Alexei Lvovich Semenov\index{Semenov, A.L.}
(1950--- ) and A.Kh. Shen',\index{Shen', A.Kh.}
{\em Russ.\ Math.\ Surveys}, 45:1(1990), 121-189].

\section{Randomness as Membership of All Large Majorities}
\label{sect.expect}
For a better understanding of the problem
revealed by Ville, as in Section~\ref{sect.unpredictability},
 and its subsequent
solution by P. Martin-L\"{o}f\index{Martin-L\"of, P.} in 1966,
we look at some aspects of the methodology
of probability theory.

\subsection{Typicality}
Consider the sample space of all one-way infinite binary sequences
generated by
fair coin tosses.
Intuitively, we call a sequence `random' iff it is `typical'.
\index{sequence!typical}
It is not `typical', say `special', if it has
a particular distinguishing property.
An example of such a property is that an infinite sequence
contains only finitely many ones. There are infinitely
many such sequences. But the probability
that such a sequence occurs as the outcome of fair coin tosses
is zero. `Typical' infinite sequences will have the converse
property, namely, they contain infinitely many ones.

In fact, one would like to say that `typical' infinite
sequences will have {\em all converse properties} of the properties
which can be enjoyed by `special' infinite sequences.
That is, such a sequence should belong to {\em all large majorities}.
This can be formalized as follows.

Suppose that a single particular property,
such as containing infinitely many occurrences of ones (or zeros),
the Law of Large Numbers\index{Law!of Large Numbers},
or the Law of
the Iterated Logarithm,\index{Law!Iterated Logarithm} has been shown to
have probability one, then one calls this
a {\it Law of Randomness}. Each sequence satisfying this property
belongs to a large majority, namely the set of all sequences satisfying
the property which has measure one by our assumption.

Now we call an infinite sequence is `typical' or `random' if it
belongs to {\em all majorities of measure one}, that is, it satisfies
all Laws of Randomness. 
In other words, each {\em single} individual `random' 
infinite sequence 
posesses all properties which hold with probability
one for the ensemble of {\em all} infinite sequences. This
is the substance of so-called pseudo-randomness tests.
For example, to test whether the sequence of digits
corresponding to the decimal expansion of $\pi = 3.1415 \ldots$
is random one tests whether the initial segment satisfies
some properties which hold with probability one
for the ensemble of all sequences.  

\begin{example}
\rm
One such
property is `normality' of a sequence as defined earlier.
Around 1909, Emile Borel called an
infinite sequence of decimal digits {\em normal} in the scale of ten
if, for each $k$, the frequency of occurrences
(possibly overlapping)
of each block $y$ of length $k \geq 1$ in the initial segment of length $n$
goes to limit $10^{-k}$ for $n$ grows unbounded,
\cite{Kn81}.
It is known that normality is not sufficient for randomness,
since Champernowne's sequence 
\[123456789101112 \ldots\]
is normal in the scale of ten. On the other hand, 
it is universally agreed that a random infinite sequence
must be normal. (If not, then some blocks occur
more frequent than others, which can be used to
obtain better than fair odds for prediction.)

For a particular binary sequence $\omega = \omega_1 \omega_2 \ldots$ let
$f_n = \omega_1 +  \omega_2 + \cdots + \omega_n$.
Of course, we cannot effectively test an infinite sequence.
Therefore, a so-called pseudo-randomness test examines
increasingly long initial segments of the individual sequence 
under consideration.

We can define a pseudo randomness test for the normality property with $k=1$
to test a candidate infinite sequence for 
increasing $n$ whether the deviations from one half
0's and 1's become too large. For example, by checking for each
successive $n$ whether
\[ | f_n - \frac{n}{2}| > (1+\epsilon) \sqrt { \frac{n \log \log n}{2}}\]
for a fixed constant $\epsilon > 0$.
(The Law of the Iterated Logarithm states that this inequality
should not hold for infinitely many $n$). If within $n$ trials
in this process we find that the inequality holds
$k$ times, then we assume the original infinite sequence to be random
with confidence at most, say, $\sum_{i=1}^n 1/2^i - \sum_{i=1}^k 1/2^i$.
(The sequence is random if the confidence is greater than zero for 
$n$ goes to infinity, and not random otherwise.)

Clearly, the number of pseudo-randomness tests we can devise
is infinite. Namely, just for the normality property alone there
is a similar pseudo-randomness test
for each $k \geq 1$.
\end{example}

But now we are in trouble. 
The naive execution of the above ideas
in classical mathematics is infeasible.
Each
individual infinite sequence induces its very own pseudo-randomness test
which tests whether a candidate infinite sequence is in fact
that individual sequence.  Each infinite sequence forms a singleton
set in the sample space of all infinite sequences. 
{\em All} complements of singleton sets
in the sample space
have probability one. 
The intersection of all complements of singleton sets is
clearly empty. Therefore, 
the intersection
of all sets of probability one is empty. Thus,
there are no random infinite sequences!

Let us give a concrete example.
Consider as sample space $S$ the set
of all one-way infinite binary sequences.
The cylinder $\Gamma_x =   \{  \omega : \omega = x  \ldots   \}  $
consists of all infinite binary sequences starting with
the finite binary sequence $x$. For instance,
$\Gamma_{\epsilon} = S$, where $\epsilon$ denotes the {\em empty sequence},
that is, the sequence with zero elements. The uniform distribution $\lambda$
on the sample space is defined by
$\lambda ( \Gamma_x ) = 2^{-l(x)}$.
That is, the probability of an
infinite binary sequence $\omega$
starting with a finite initial segment $x$ is $2^{-l(x)}$.
In probability theory it is
general practice that if a certain property,
such as the Law of Large Numbers\index{Law!of Large Numbers}, or the Law of
the Iterated Logarithm\index{Law!Iterated Logarithm}, has been shown to
have probability one, then one calls this
a %
\it Law of Randomness\index{Law!of Randomness|bold}%
\rm . For example,
in our sample space the Law of Large Numbers
says that $\lim_{{n}   \rightarrow   \infty}  ( \omega_1 +  \cdots   + \omega_n
)/n = 1/2$.
If $A$ is the set of elements of $S$ which satisfy the Law of
Large Numbers, then it can be shown that $\lambda (A) = 1$.

Generalizing this idea for $S$ with measure $\mu$, one may identify %
\it any %
\rm set $B   \subseteq   S$,
such that $\mu (B) = 1$, with the Law of Randomness, namely,
`to be an element of $B$'. Elements of $S$
which are do not satisfy the law `to be an element of
$B$' form a set of measure zero, a
\it
null set\index{null set|bold}.
\rm
It is natural to call an element
of the sample space `random' if it satisfies all laws
of randomness. Now we are in trouble.
For each element $\omega   \in   S$,
the set $B_{\omega} = S -   \{  \omega  \}  $
forms a law of randomness.
But the intersection
of all these sets $B_{\omega}$ of probability one is empty. Thus,
no sequence would be random, if we require
that all laws of randomness that `exist'
are satisfied by
a random sequence.

\subsection{Randomness in Martin-L\"of's Sense}
The Swedish mathematician
Per Martin-L\"of\index{Martin-L\"of, P.},
visiting Kolmogorov in Moscow during 1964-1965, 
investigated the complexity oscillations of infinite sequences
and proposed a definition of infinite random sequences
which is based on constructive measure theory
using ideas related to Kolmogorov
complexity, 
[%
\it Inform.\ Contr.%
\rm , 9(1966), 602-619;
\it Z.\ Wahrsch.\ Verw.\ Geb.,
\rm 19(1971), 225-230].
This way he succeeded
in defining random infinite sequences in a manner which is free
of above difficulties. 

It turns out that a
constructive viewpoint enables us to carry
out this program mathematically without such pitfalls.
In practice, all laws
that are proved in probability theory to hold
with probability one are effective in the formal sense
of effective computability due to A.M. Turing.
A straightforward formalization of this viewpoint is to require a
law of probability to be partial recursive in the sense that we can
effectively test whether it is violated.
This suggests that the set of random infinite sequences
should not be defined as the
intersection of all sets of measure one,
but as the intersection of all
sets of measure one with a recursively enumerable complement.
The latter intersection
is again a set of measure one with a recursively enumerable complement.
Hence, there is a single effective law of randomness
which can be stated as the property `to satisfy all effective
laws of randomness', and the infinite sequences
have this property
with probability one.

(There is a related development in set theory and recursion theory,
namely, the notion of `generic object' in the context of `forcing'.
For example, an equivalent definition of genericity in arithmetic
is being a member of the intersection of all arithmetical sets
of measure 1. There is a notion, called `1-genericity',
which calls for the intersection of all
recursively enumerable sets of measure 1.
This is obviously related to the approach of Martin-L\"of,
and prior to it.
Forcing was introduced by Paul Cohen (1934--- )\index{Cohen, P.}
in 1963 to show the independence of the continuum hypothesis,
and using sets of positive measure as forcing conditions is
due to Robert M. Solovay\index{Solovay, R.M.} soon afterwards.)

The natural formalization is to identify the
effective test with a partial recursive function.
This suggests that one ought to consider not the
intersection of all sets of measure one,
but only the intersection of all
sets of measure one with recursively enumerable
complements.
(Such a complement set is expressed as the
union of a recursively enumerable set of cylinders)\index{cylinder}.
It turns out that this intersection
has again measure one. Hence, almost all infinite sequences
satisfy all effective Laws of Randomness
with probability one.
This notion of infinite random sequences turns out to be related to
infinite sequences of which all finite initial segments
have high Kolmogorov complexity.
\begin{comment}
The notion of randomness satisfied by
both the Mises-Wald-Church\index{sequence!Mises-Wald-Church random}
collectives and the Martin-L\"of random infinite sequences
\index{sequence!random}
is roughly that
\it effective tests
\rm cannot
detect regularity. This does not mean that
a sequence may not exhibit regularities which
cannot be effectively tested. Collectives generated by Nature,
\index{collective}
as postulated by von Mises, may very well always
\index{Mises, R. von}
satisfy stricter criteria of randomness.
Why should collectives generated by quantum mechanic
phenomena care about
mathematical notions
of computability? Again,
satisfaction of all effectively testable prerequisites
for randomness is some form of regularity.
Maybe nature is more lawless than adhering strictly
to regularities imposed by the statistics of randomness.
\end{comment}
Until now the discussion has centered on
infinite random sequences where the randomness is defined in terms
of limits of relative frequencies. However,
\begin{comment}
``The frequency concept based on the notion of
\it limiting
frequency %
\rm as the number of trials increases to infinity,
does not contribute anything to substantiate the application
of the results of probability theory to real practical
problems where we always have to deal with a finite number
of trials.'' 
[A.N. Kolmogorov,
On tables of random numbers,
$Sankhy \bar a$, Series A, 25(1963), 369-376. Page 369.]
\index{Kolmogorov, A.N.}
\end{comment}
The practical objection against both the relevance
of considering infinite sequences of trials
and the existence of a relative frequency limit
is concisely put in John Maynard  Keynes' (1883---1946)\index{Keynes, J.M.} 
famous phrase
``in the long run we shall all be dead.''
It seems more appealing to try to
define randomness for finite strings first, and only
then define random infinite strings in terms of randomness of initial
segments.

The approach of von Mises to define randomness of infinite sequences
\index{Mises, R. von}
in terms of %
\it unpredictability %
\index{prediction error}
\rm of continuations
of finite initial sequences under certain laws (like
recursive functions) did not lead to satisfying results.
The Martin-L\"of approach does lead to satisfying results,
and is to a great extent equivalent with the
Kolmogorov complexity approach.
Although certainly inspired by the random sequence debate, the
introduction of Kolmogorov complexity marks a definite
shift of point of departure. Namely, to
define randomness of sequences by the fact that no
program from which an initial segment of the sequence can
be computed is significantly shorter\label{compl.initial.segm}
than the initial segment itself, rather than that no program can
predict the next elements of the sequence.
Thus, we change the focus from the `unpredictability' criterion
to the `incompressibility' criterion, and since this will
turn out to be equivalent with Martin-L\"of's approach,
the `incompressibility' criterion is both necessary and
sufficient.

\subsection{Random Finite Sequences}\label{sect.random.finite}
Finite sequences which cannot be effectively described in
a significant shorter description than their literal representation
are called random.
Our aim is to characterize random infinite
sequences as sequences of which all initial finite segments
are random in this sense.
Martin-L\"of's related approach
characterizes random infinite
sequences as sequences of which all initial finite segments
pass all effective randomness tests.
\begin{comment}
Initially, before the idea of complexity, Kolmogorov
\index{Kolmogorov, A.N.}
proposed a
close analogy to von Mises'
\index{Mises, R. von}
notions in the finite domain.
Consider a generalization
of place-selection rules
\index{place-selection rule!according to Kolmogorov}
insofar as the selection of
$a_i$ can depend on $a_j$ with $j > i$ [A.N. Kolmogorov,
$Sankhy \bar a$, Series A, 25(1963), 369-376].
Let $\Phi$ be a finite set of such generalized place-selection rules.
Kolmogorov suggested that an arbitrary finite binary
sequence $a$ of length $n \geq m$ can be called $(m, \epsilon )$-random
with respect to $\Phi$, if there exists some $p$ such that
the relative frequency of the 1's in the subsequences
$a_{i_1 } \ldots a_{i_r}$ with $r \geq m$, selected
by applying some $\phi$ in $\Phi$ to $a$, all
lie within $\epsilon $ of $p$. (We discard $\phi$ that
yield subsequences shorter than $m$.) Stated differently,
the relative frequency in this finite subsequence is
approximately (to within $\epsilon$) invariant
under any of the methods of subsequence selection that
yield subsequences of length at least $m$. Kolmogorov has
\index{Kolmogorov, A.N.}
shown that if the cardinality of $\Phi$ satisfies:
$$
d( \Phi )  \leq {1 \over 2 }e^{{2m} \epsilon^2 (1 - \epsilon )} ,
$$
then, for any $p$ and any $n \geq m$ there is some sequence $a$
of length $n$ which is $(m, \epsilon )$-random with respect
to $\Phi$.
\end{comment}

Let us borrow some ideas from statistics.
We are given a certain sample space $S$ with an associated distribution $P$.
Given an element $x$ of the sample space,
we want to test the hypothesis `$x$ is a typical outcome'.
Practically speaking, the property of being typical is
the property of belonging to any reasonable majority.
In choosing an object at random, we have confidence that
this object will fall precisely in the intersection
of all such majorities. The latter condition we identify
with $x$ being random.

To ascertain whether a given element of the sample space
belongs to a particular reasonable majority we introduce the notion of a test.
\index{test!in statistics|bold}
Generally, a test is given by a prescription which, for
every level of significance $\epsilon$, tells us for what
elements $x$ of $S$ the hypothesis `$x$ belongs to majority $M$ in $S$'
should be rejected, where $\epsilon = 1 - P(M)$.
Taking $\epsilon = 2^{{-} m}$, $m = 1, 2, \ldots $,
this amounts to saying that we have a description
of the set $V   \subseteq   {\cal N} \times S$ of nested %
\it critical regions %
\index{test!critical region of|bold}
\rm
\begin{eqnarray*}
 V_m & = &  \{  x: (m, x)   \in   V  \} \\
 V_m  &  \supseteq  & V_{{m} + 1} ,   \ \ \ \  m = 1, 2,  \ldots .
\end{eqnarray*}
The condition that $V_m$ be a critical region on the
{\it significance level}\index{test!significance level of}
$\epsilon  = 2^{-m}$ amounts to requiring, for all $n$
$$
\sum_x    \{  P(x): l(x) = n, x   \in   V_m   \}    \leq \epsilon .
$$
The complement of a critical region $V_m$ is called the
$(1- \epsilon )$ {\it confidence interval}.
\index{test!confidence interval of}
If $x \in V_m$, then
the hypothesis `$x$ belongs to majority $M$', and therefore
the stronger hypothesis `$x$ is random', is rejected with
\index{test!testing for randomness}
significance level $\epsilon $.
We can say that $x$ fails the test at the level of
critical region $V_m$.

\begin{example}\label{example.MLtest}
\rm
A string $x_1 x_2  \ldots x_n$
with many initial zeros is not very
random. We can test this aspect as follows.
The special test $V$ has
critical regions $V_1 , V_2 , \ldots $.
Consider $x = 0.x_1 x_2  \ldots x_n$ as a rational number,
and each critical region as a half-open interval
$V_m = [0, 2^{-m}  )$ in $[0, 1)$, $m = 1, 2, \ldots $.
Then the subsequent
critical regions test the hypothesis `$x$ is random' by considering
the subsequent digits in the binary expansion of $x$.
We reject the hypothesis on the
significance level $\epsilon = 2^{-m}$
provided $x_1 = x_2 = \cdots = x_m = 0$,
Another test for randomness of finite binary
strings rejects when the relative frequency
of ones differs too much from $1/2$.
This particular test can be implemented by
rejecting the hypothesis of randomness of
$x=x_1 x_2 \ldots x_n$ at level $\epsilon  = 2^{-m}$ provided
$|2 f_n - n|   >   g(n, m)$,
where $f_n = \sum_{i=1}^n x_i$,
and $g(n, m)$ is the least number determined by the requirement
that the number of binary strings $x$
of length $n$ for which this inequality holds
is at most $2^{n-m}$.
\end{example}

In practice, statistical tests are %
\it effective %
\rm prescriptions
such that we can compute, at each level of significance,
for what strings the associated hypothesis should be rejected.
It would be hard to imagine what use
it would be in statistics to have tests
that are not effective in the sense of computability
theory.

\begin{definition}
\rm
\label{test}
Let $P$ be a recursive
probability distribution
on the sample space ${\cal N}$.
A total\index{test|bold}
function $\delta : {\cal N}   \rightarrow   {\cal N}$ is
a $P$-{\em test} (Martin-L\"of
test\index{test!Martin-L\"of}
for randomness)\index{test!$P$-} if:
\begin{enumerate}
\item
$\delta$ is enumerable
(the set
$V =   \{  (m, x): \delta (x)   \geq   m  \}  $ %
\rm is recursively
enumerable); and
\item
$\sum   \{  P(x): \delta (x)   \geq   m ,$
$l(x) =$
$n   \}    \leq $
$2^{-m}$, for all $n$.
\end{enumerate}
\end{definition}

The critical regions associated with
the common statistical tests
are present in the form
of the sequence
$V_1   \supseteq   V_2   \supseteq    \cdots $, where
$V_m =   \{  x: \delta (x)   \geq   m  \}  $%
\rm , for $m   \geq   1$.
Nesting is assured since $\delta (x)   \geq   m+1$ implies
$\delta (x)   \geq   m$. Each set $V_m$ is recursively enumerable
because of Item 1.

A particularly important case is $P$ is
the uniform distribution\index{distribution!uniform}, defined by
$L(x) = 2^{-2l(x)}$. The restriction of $L$ to strings of length
$n$ is defined by $L_n (x)=2^{-n}$ for $l(x)=n$ and 0 otherwise.
(By definition, $L_n (x)=L(x|l(x)=n)$.)
Then, Item 2 can be rewritten as $\sum_{x \in V_m} L_n (x) \leq 2^{-m}$
which is the same as
$$
d(  \{  x: l(x) = n, \: x   \in   V_m    \}  )  \leq 2^{n-m}.
$$
In this case we often speak simply of a %
\it test\index{test}%
\rm ,
with the uniform distribution $L$ understood.

\begin{comment}
In statistical
tests
membership of $(m,x)$
in $V$ can usually be determined in time polynomial in $l(m) + l(x)$.
\end{comment}

\begin{example}\label{example.test.oddones}
\rm
The previous test examples can be rephrased
in terms of Martin-L\"of tests.
Let us try a more subtle example. A real number such that
all bits in odd positions in its binary representation are 1's
is not random with respect to the uniform distribution.
To show this we need a test which detects
sequences of the form $ x =  1 x_2 1 x_4 1 x_6 1 x_8   \ldots $.
Define a test $\delta$ by
$$
\delta (x)  = \max   \{   i:  x_1 = x_3 =  \cdots  = x_{2i-1} = 1  \}  ,
$$
and $\delta (x) = 0$ if $x_1 = 0$.
For example: $\delta (01111) = 0$; $\delta (10011) = 1$; $\delta (11011) = 1$;
$\delta (10100) = 2$; $\delta (11111) = 3$.
To show that $\delta$ is a test we have to show that $\delta$
satisfies the definition of a test. Clearly, $\delta$ is
enumerable (even recursive).
If $\delta (x)   \geq   m$ where $l(x) = n   \geq   2m$,
then there are $2^{m-1}$ possibilities for the $(2m - 1)$-length
prefix of $x$, and $2^{n-(2m-1)}$ possibilities for the remainder
of $x$. Therefore,
$d  \{  x: \delta (x)   \geq   m,  \: l(x) = n  \}   \leq 2^{n-m}$.
\end{example}

\begin{definition}\label{def.martinloeftest}
\rm
A universal Martin-L\"of test for randomness with respect
to distribution $P$, a %
\it universal P-test %
\index{test!universal|bold}
\rm for short,
is a test $\delta_0 (\cdot |P)$ %
\rm  such that for
each $P$-test $\delta$, there is a constant $c$,
such that for all $x$, we have $\delta_0 (x|P)   \geq   \delta (x) - c$.
\end{definition}

We started out with the objective to establish in what sense
incompressible strings may be called random.
Kolmogorov considered a notion of
{\em randomness deficiency} $\delta (x| A) = l(d(A))-C(x|A)$ of a string $x$
relative to  a finite set $A$.
With $A$ the set of strings of length $n$
and $x \in A$ we find
$\delta (x| A) = \delta (x|n) = n - C(x| n)$.

\begin{theorem}
\label{M1}
The function
$f(x) = n - C(x| n)-1$
is a universal $L$-test with $L$ the uniform distribution
over $\{0,1\}^*$ and $n=l(x)$.
\end{theorem}

\subsection{Random Infinite Sequences}

\index{sequence!random|(}
\label{sect.random.infinite}
Consider the question of how $C$ behaves
in terms of increasingly long initial segments of a fixed
infinite binary sequence (or string) $\omega$. For instance,
is it monotone in the sense that $C(\omega_{1:m} )  \leq C(\omega_{1:n} )$,
or $C(\omega_{1:m} | m)  \leq C(\omega_{1:n} | n)$,
for all infinite binary sequences $\omega$ and all $m  \leq n$?
We have already seen
that the answer is negative in both cases.
A similar effect arises when we try to use Kolmogorov complexity
to solve the problem of finding
a proper definition of random infinite sequences
({\em collectives})
\index{collective}
according to the task already set by von Mises\index{Mises, R. von} in 1919,
Section~\ref{sect.unpredictability}.

It is seductive
to call an infinite binary sequence $\omega$ random if
there is a constant $c$ such that, for all $n$, the $n$-length
prefix $\omega_{1:n}$ has $C(\omega_{1:n} )   \geq   n-c$. However,
such sequences
do not exist.

As in Section~\ref{sect.random.finite}, we define
a test for randomness. However, this time the test
will not be defined on the entire sequence (which is
impossible for an effective test and an infinite sequence),
but for each finite binary string. The value of
the test for an infinite sequence is then defined
as the maximum of the values of the test on
all prefixes. Since this suggests an effective
process of sequential approximations, we
call it a sequential test. Below, we need to use
notions of continuous sample spaces and measures.

\begin{definition}\label{def.sequential.test}
\rm
Let $\mu$ be a recursive probability measure on
the sample space
$  \{  0,1  \} ^{\infty} $.
A total function
$\delta :   \{  0, 1  \} ^{\infty}   \rightarrow
 {\cal N}   \cup     \{  \infty  \}  $
is a\index{test!sequential Martin-L\"of}
{\it sequential $\mu$-test}
(sequential Martin-L\"of\index{test!sequential $\mu$-}
$\mu$-test for randomness)\index{test!sequential|bold}
if:
\begin{enumerate}
\item
$\delta ( \omega ) = \sup_{n \in {\cal N}}\{  \gamma  ( \omega_{1:n} ) \}  $,
where $\gamma  : {\cal N}   \rightarrow   {\cal N}$
is a total enumerable function
($V = \{  (m, y): \gamma  (y) \geq m  \}$
is a recursively enumerable set); and
\item
$\mu  \{ \omega : \delta ( \omega ) \geq m \} \leq 2^{-m}$, for each $m   \geq
 0$.
\end{enumerate}

If $\mu$ is the uniform measure $\lambda$, then we often use simply
\it sequential
test\index{test!sequential for uniform distribution}%
\rm .
\end{definition}

\begin{comment}
We can require $\gamma$ to be a recursive function without changing
the notion of a sequential $\mu$-test. By definition, for each enumerable
function $\gamma$ there exists a recursive function
$\phi$ such that $\phi (x,k)$ nondecreasing in $k$ such
that $\lim_{k \rightarrow \infty } \phi (x,k) = \gamma (x)$.
Define a recursive function $\gamma '$ by
$\gamma ' (\omega_{1:n})= \phi (\omega_{1:m},k)$ with
$\langle m,k \rangle =n$. Then,
$\sup_{n \in {\cal N}} \{ \gamma ' (\omega_{1:n} ) \} =
\sup_{n \in {\cal N}} \{ \gamma (\omega_{1:n} ) \}$.
\end{comment}

\begin{example}\label{ex.seq.test.1}
\rm
Consider $\{0,1\}^{\infty}$ with the uniform measure
$\lambda (x) = 2^{-l(x)}$. An example of a sequential $\lambda$-test is
to test whether there are 1's in even positions
of $\omega \in \{0,1\}^{\infty}$. Let
\[ \gamma (\omega_{1:n} ) =  \left\{ \begin{array}{ll}
n & \mbox{if $\sum_{i=1}^{n/2} \omega_{2i} = 0$,} \\
0 & \mbox{otherwise.}
\end{array}
\right. \]
The number of $x$'s of length $n$ such that
$\gamma (x) \geq m$ is at most $2^{n/2}$ for any $m \geq 1$. Therefore,
$\lambda \{ \omega: \delta (\omega ) \geq m \} = 0 \leq 2^{-m}$ for
$m>0$. For $m=0$,
$\lambda \{ \omega: \delta (\omega ) \geq m \} \leq 2^{-m}$ holds
trivially. A sequence $\omega$ is random with respect to this test
if $\delta (\omega ) < \infty$. Thus, a sequence $\zeta$
with 0's in all even locations will have $\delta (\zeta ) = \infty$,
and it will fail the test and hence $\zeta$ is not random
with respect to this test. Notice that this is not a very strong test
of randomness. For example, a sequence $\eta = 010^{\infty}$ will
pass $\delta$ and be considered random with respect to this
test. This test only filters out some nonrandom sequences with
all 0's at the even locations and cannot detect other
kinds of regularities.
\end{example}

We continue the general theory of sequential testing.
If $\delta (\omega) = \infty$, then we say that $\omega$ %
\it fails %
\rm $\delta$, or that $\delta$ %
\it rejects %
\rm $\omega$.
Otherwise, $\omega$ %
\it passes %
\rm $\delta$.
By definition, the set of $\omega$'s that is rejected by $\delta$
has $\mu$-measure zero, and, conversely, the set of $\omega$'s that
pass $\delta$ has $\mu$-measure one.

Suppose for a test $\delta$ holds $\delta ( \omega) = m$.
Then there is a prefix $y$ of $ \omega$, with $l(y)$ minimal,
such that $\gamma  (y) = m$ for the $\gamma $ used to define $\delta$.
Then, obviously, %
\it each %
\rm infinite sequence $\zeta$ that starts with $y$
has $\delta ( \zeta)   \geq   m$. The set of such $\zeta$ is
$\Gamma_y =$
$  \{  \zeta : \zeta = y \rho ,   \rho \in \{  0, 1  \} ^{\infty}  \}  $,
the cylinder\index{cylinder} generated by $y$. Geometrically speaking,
$\Gamma_y$ can be viewed as the set of all real numbers $0.y \ldots $
corresponding to the half-open %
\it interval %
\rm $I_y = [0.y, 0.y+2^{-l(y)} )$.
(For the uniform measure, $\lambda( \Gamma_y ) = 2^{-l(y)}$,
the common length of $I_y$.)

In terms of common statistical tests, the critical regions are
formed by the nested sequence
$V_1    \supseteq   V_2    \supseteq    \cdots $,
where $V_m$ is defined as
$V_m  =   \{  \omega : \delta ( \omega)   \geq   m  \}  $, for $m   \geq   1$.
We can formulate the definition of $V_m$ as
$$
V_m = \bigcup   \{  \Gamma_y : (m, y)   \in   V  \}  .
$$
In geometric terms, $V_m$ is the union of
a set of subintervals of $[0, 1)$.
Since $V$ is recursively enumerable, so is
the set
of intervals whose union is $V_m$.
For each critical section we
have $\mu(V_m )  \leq 2^{-m}$ (in the measure
we count overlapping intervals only once).

Now we can reformulate the notion of passing a sequential test $\delta$
with associated set $V$:
$$
\delta ( \omega )   <   \infty {\rm \ iff\ } \omega  \not\in
\bigcap_{{m} = 1}^{\infty} V_m .
$$

\begin{definition}\label{def.rand.inf.seq}
\rm
\index{sequence!random|bold}
Let ${\bf V}$ be the set of all sequential $\mu$-tests.
An infinite binary sequence $ \omega$,
or the binary represented real number $0. \omega$,
is called %
\it $\mu$-random %
\rm if it passes %
\it all %
\rm sequential $\mu$-tests:
$$
\omega \not\in \bigcup_{{V}   \in   {\bf V}} \hspace{0.5em}
  \bigcap_{{m} = 1}^{\infty}  V_m .
$$
\end{definition}

For each sequential $\mu$-test $V$,
we have $\mu( \bigcap_{{m} = 1}^{\infty}  V_m ) = 0$,
by Definition~\ref{def.sequential.test}. We
call $\bigcap_{{m} = 1}^{\infty}  V_m$ a %
\it constructive $\mu$-null set%
\index{null set!constructive|bold}
\index{null set!$\mu$-}
\rm .
Since there are only
countably infinitely many sequential $\mu$-tests $V$, it follows from
standard measure theory that
$$
\mu \left( \bigcup_{{V} \in {\bf V}} \hspace{0.5em}
 \bigcap_{{m} = 1}^{\infty} V_m \right) =0,
$$
and we call the set
$U  = \bigcup_{{V}   \in   {\bf V}} \bigcap_{{m} = 1}^{\infty}  V_m$
the %
\it maximal constructive $\mu$-null set
\rm .
\\

Similar to Section~\ref{sect.random.finite}, we construct
an enumerable function $\delta_0 ( \omega |\mu )$, the
universal sequential $\mu$-test\index{test!universal sequential}
which incorporates (majorizes)
all sequential $\mu$-tests $\delta_1 , \delta_2 , \ldots $,
and that corresponds to $U$.

\begin{definition}\label{def.universalsequentialtest}
\rm
A %
{\em universal sequential $\mu$-test}\index{test!universal sequential|bold}
\rm $f$ is a
sequential $\mu$-test such that for each sequential
$\mu$-test $\delta_i$ there is a constant $c   \geq   0$ and
for all $\omega   \in     \{  0, 1  \} ^{\infty}$, we have
$f ( \omega )   \geq   \delta_i ( \omega ) - c$.
\end{definition}

\begin{theorem}
There is a universal sequential $\mu$-test {\rm (}denoted
as $\delta_0 (\cdot |\mu)${\rm )}.
\end{theorem}
\begin{definition}\label{def.random.sequence}
\rm
Let the sample space $\{0,1\}^{\infty}$ be distributed
according to $\mu$, and
let $\delta_0(\cdot |\mu)$ be a universal sequential $\mu$-test.
An infinite binary sequence $\omega$ is %
{\em $\mu$-random
in the sense of Martin-L\"of},
if $\delta_0 ( \omega |\mu )   <   \infty$.
We call such a sequence simply {\em random},
where both $\mu$ and Martin-L\"of are understood.
(This is particularly interesting for $\mu$ is
the uniform measure.)
\end{definition}

Note that %
\it this definition does not depend
on the choice of the particular universal sequential $\mu$-test %
\rm with respect to which the level is defined.
Hence, the line between random and nonrandom
infinite sequences is drawn
sharply without dependence on a reference $\mu$-test.
It is easy to see that the set of infinite sequences, which are
not random in the sense of Martin-L\"of, forms
precisely the maximal constructive $\mu$-null set of
\index{null set!constructive}
$\mu$-measure zero we have constructed above. Therefore,
\begin{theorem}
\label{M2}
Let $\mu$ be a recursive measure.
The set of $\mu$-random infinite binary sequences has $\mu$-measure one.
\end{theorem}

We say that the universal sequential $\mu$-test $\delta_0(\cdot |\mu)$
rejects an infinite sequence with probability zero,
and we conclude that a randomly selected infinite sequence passes all
effectively testable laws of randomness with probability one.

\begin{comment}
The main question remaining is the following.
Let $\lambda$ be the uniform measure.
Can we formulate a universal
sequential $\lambda$-test\index{test!universal sequential
for the uniform measure}
in terms of
complexity?
In Theorem~\ref{M1} the universal
(nonsequential) test is expressed that way.
The most obvious candidate for the
universal sequential test would be
$f ( \omega ) = \sup_{n \in {\cal N}} \{  n - C( \omega_{1:n} ) \}  $,
but it is improper. To see this,
it is simplest to notice that $f( \omega )$ would declare
\it all %
\rm infinite $ \omega $ to be nonrandom since
$f( \omega ) = \infty$, for all $ \omega $, because
of the complexity oscillations we have discussed above.
The same would be the case for
$f ( \omega ) = \sup_{n \in {\cal N}} \{  n - C( \omega _{1:n} | n ) \}  $,
by about the same proof. It is difficult to
express a universal sequential test precisely in terms of $C$-complexity.
Yet it is easy to separate the random infinite sequences from the nonrandom ones
in terms of `prefix' complexity, see below.
\end{comment}
The idea that random infinite sequences are those sequences
such that the complexity of the initial $n$-length segments
is at most a fixed additive constant below $n$, for all $n$, is one
of the first-rate ideas in the area of Kolmogorov
complexity. In fact, this was one of the motivations for
Kolmogorov to invent Kolmogorov complexity in the first
place. We have seen
in Section~\ref{sect.random.infinite} that this does not
work for the plain Kolmogorov complexity $C(\cdot)$,
due to the complexity oscillations. The next result is
important, and is a culmination of the theory.
For prefix complexity $K(\cdot)$ it is indeed the case
that random sequences are those sequences for which the complexity
of each initial segment is at least its length.
\begin{theorem}
\label{K3}
An infinite binary sequence $ \omega$
is random with respect to the uniform measure
iff
there is
a constant $c$ such that 
$K( \omega_{1:n} )   \geq   n - c$,
for all $n$.
\end{theorem}
Although the following discussion is a bit beyond the scope
of this article, it is important to understand the raised issues.
There are different `families' of tests which characterize precisely
the same class of random infinite sequences. The sequential
tests are but one type. The test in Theorem~\ref{K3} is an example of
an integral test (a universal one) with respect
to the uniform measure as defined.
The introduction of different
families of tests like  martingale tests and integral tests
requires special machinery. It is of the utmost importance
to realize that an infinite sequence is random only in the sense
of being random with respect to a given distribution. For
example, if we have the sample space $\{0,1\}^{\infty}$
and a probability measure $\mu$ such that $\Gamma_x$
has $\mu$-measure one for all sequences $x=00 \ldots 0$
and measure zero for all $\Gamma_y$ with $y \neq 00 \ldots 0$,
then the only $\mu$-random infinite sequence is $\omega = 00 \ldots$.

Thus, in many applications we are not really interested in randomness
with respect to the uniform distribution, but 
in randomness with respect to a given recursive distribution $\mu$.
It is therefore important to have explicit expressions for
$\mu$-randomness. In particular, we have
\begin{corollary}\label{K3corgeneral}
\rm
Let $\mu$ be a recursive measure. The function
\[\rho_0 ( \omega |\mu )= \sup_{\omega \in \Gamma_x}
 \{ - K(x| \mu) - \log \mu ( x) \} \]
is a universal integral $\mu$-test.
\end{corollary}

\begin{example}
\rm
With respect to the special case
of the uniform distribution $\lambda$, Theorem~\ref{K3}
sets $\rho_0 (\omega |\lambda)= \sup_{n \in {\cal N}} \{n-K(\omega_{1:n})\}$ up
to a
constant additional term. This is the expression
we found already in  Theorem~\ref{K3}.
\end{example}

Such theories of exact expressions
for tests which separate the random infinite sequences
from the nonrandom ones with respect to any recursive measure
are presented in \cite{LV90}, together with applications
in inductive reasoning (the minimum description length
inference in statistics) and physics and computation.
In the latter, one expresses a variation of coarse-grained Boltzmann entropy 
of an {\em individual}
micro state of the system in terms of a generalization of Kolmogorov
complexity called `algorithmic entropy' rather than conventional
coarse-grained Boltzmann entropy in statistical mechanics
which expresses the uncertainty of the micro state
when it is constrained only by the macro state comprising
an ensemble of micro states.

\begin{example}\label{ex.halting.probability}
\rm
\index{halting probability|(}
It is impossible to {\em construct} an infinite
random sequence by algorithmic means. But,
using the reference prefix machine
$U$, we can {\em define} a particular
random infinite binary sequence
in a more or less natural way.
The {\em halting probability} is the real number $\Omega$
defined by
\index{halting probability|bold}\index{$\Omega$: halting probability}
$$
\Omega = \sum_{{U(p)}   <   \infty}  2^{{-} l(p)},
$$
the sum taken over all inputs $p$ for which the reference
machine $U$ halts.
Since $U$ halts
for some $p$, we have $\Omega   >   0$.
Because $U$ is a prefix machine, the set
of its programs forms a prefix-code and
by Kraft's Inequality, see \cite{LV90}, we find
$\Omega  \leq 1$. Actually, $\Omega   <   1$ since
$U$ does not always halt.

We call $\Omega$ the halting probability because it is
the probability that $U$ halts if its
program is provided by a sequence of fair coin flips.
G.J. Chaitin observed that the number $\Omega$ has interesting properties.
In the first place, the binary representation of
the real number $\Omega$ encodes
the halting problem\index{halting problem}
very compactly.
Denote the initial $n$-length
segment of $\Omega$ after the decimal dot by $\Omega_{1:n}$.
If $\Omega$ were a terminating binary rational number, then
we use the representation with infinitely many zeros,
so that $\Omega < \Omega_{1:n} + 2^{-n}$.

\begin{claim} \label{claim.omega1}
\rm
Let $p$ be a binary string of length at most $n$.
Given $\Omega_{1:n}$, it is decidable whether the
reference prefix machine $U$ halts on input $p$.
\end{claim}

\begin{proof}
Clearly,
\begin{equation}
\Omega_{1:n}    \leq    \Omega   < \Omega_{1:n}  +  2^{-n}.
\label{omega1}
\end{equation}
Dovetail the computations of $U$ on all inputs as follows.
The first phase consists of $U$ executing one step
of the computation on the first input. In the second phase,
$U$ executes the second step of the computation on the first input and
the first step of the computation of the second input.
Phase $i$
consists of $U$ executing the $j$th step of the computation
on the $k$th input, for all $j$ and $k$ such that $j + k = i$.
We start with an approximation $\Omega '  := 0$. Execute
phases 1,2, \ldots . Whenever any computation of $U$ on some input $p$
terminates, we improve our approximation of $\Omega$ by
executing
$$
\Omega ' := \Omega ' + 2^{-l(p)} .
$$
This process eventually
yields an approximation $\Omega '$ of $\Omega$, such that
$\Omega '   \geq   \Omega_{1:n}$.
If $p$ is not among the halted programs
which contributed to $\Omega '$, then $p$ will never halt.
With a new $p$ halting we add
a contribution of $2^{-l(p)} \geq 2^{-n}$ to the approximation
of $\Omega$, contradicting Equation~\ref{omega1} by
\[ \Omega \geq \Omega' + 2^{-l(p)} \geq \Omega_{1:n} + 2^{-n}. \]
\end{proof}

It follows that the binary expansion of $\Omega$ is an %
\it incompressible\index{sequence!incompressible}
\rm sequence.
\begin{claim}\label{omega2}
\rm
There is a constant $c$ such that
$K( \Omega_{1:n} )   \geq   n-c$ for all $n$.
\end{claim}

\begin{comment}
\rm
That is, $\Omega$ is a
\it particular
random real\index{number!random real}
\rm , and one that is naturally defined to boot.
That $\Omega$ is random implies that it is not computable,
and therefore %
\it transcendental\index{number!transcendental}%
\rm . Namely, if it were
computable, then $K( \Omega_{1:n} | n) = O(1)$, which contradicts
Claim~\ref{omega2}.
By the way, irrationality of $\Omega$ implies
that both inequalities in Equation~\ref{omega1}
are strict.
\end{comment}

\begin{proof}
From Claim~\ref{claim.omega1} it follows that, given $\Omega_{1:n} $,
one can calculate all programs $p$ of length not greater than $n$
for which the reference prefix machine $U$ halts.
For each $x$
which is not computed by any of these halting programs
the shortest program $x^*$ has size greater than $n$,
that is, $K(x) > n $.
Hence, we can construct a recursive function $\phi$
computing such high complexity $x$'s from initial
segments of $\Omega$, such that for all $n$,
$$
K( \phi ( \Omega_{1:n} ) )    >   n .
$$
Given a description of $\phi$ in $c$ bits,
for each $n$ we can compute $\phi ( \Omega_{1:n} )$ from
$\Omega_{1:n}$, which means
$$
K( \Omega_{1:n} ) + c   \geq   n ,
$$
which was what we had to prove.
\end{proof}

\begin{corollary}
\rm
By Theorem~\ref{K3}, we find that
$\Omega$ is random with respect to the uniform measure.
\index{sequence!random}
\end{corollary}
\begin{comment}
Knowing the first 10,000 bits of $\Omega$
enables us to solve the halting of all programs of less than
10,000 bits.
This includes programs looking for counterexamples
to Fermat's Last Theorem (presumed solved affirmatively
at this time of writing),\index{Theorem!Fermat's Last}
Goldbach's Conjecture,\index{Goldbach's Conjecture}
Riemann's Hypothesis,\index{Riemann's Hypothesis} and most
other famous conjectures in mathematics which can be refuted by single
finite counterexamples. Moreover, for all axiomatic mathematical
theories which can be expressed compactly enough to be conceivably
interesting to human beings, say in less than 10,000 bits,
$\Omega_{10,000}$ can be used to decide for every statement
in the theory whether it is true, false, or independent.
Finally, knowledge of $\Omega_{1:n}$ suffices to determine
whether $K(x)  \leq n$ for each finite binary string $x$. Thus,
$\Omega$ is truly the number of Wisdom, and `can be known of,
\index{number!of Wisdom}
but not known, through human reason'
[Charles H. Bennett\index{Bennett, C.H.}
and Martin Gardner\index{Gardner, M.}, %
\it Scientific American%
\rm ,
241:11(1979), 20-34].
But even if you possess $\Omega_{1:10,000}$, you cannot
use it except by spending time of a thoroughly unrealistic
nature. (The time $t(n)$ it takes to find all halting programs
of length less than $n$ from $\Omega_{1:n}$ grows
faster than any recursive function.)
\end{comment}
\end{example}

\begin{example}\label{exam.diophantineequation}
\rm
Let us look at another example of an infinite random sequence, this time
defined in terms of
Diophantine equations. These are algebraic equations
of the form $X = 0$, where $X$ is
build up from nonnegative integer variables
and nonnegative integer constants by a finite
number of additions ($A + B$) and multiplications ($A \times B$).
The best
known examples are $x^n + y^n = z^n$, where $n = 1, 2, \ldots $.
 
\begin{comment}
Pierre de Fermat\index{Fermat, P. de} (1601-1665)
has stated that this equation
has no solution in integers $x$, $y$, and $z$ for $n$
an integer greater than 2. (For $n = 2$ there exist solutions, for instance
$3^2 + 4^2 = 5^2$.) However, he
did not supply a proof of this assertion, often
called {\em Fermat's Last Theorem}\index{Theorem!Fermat's Last|bold}.
After 350 years of withstanding concerted attacks
to come up with a proof or disproof,
the problem had become a celebrity among
unsolved mathematical problems. (At this time of writing,
October 1995, there is a serious claim that the problem has
been settled.)
Suppose we substitute all
possible values for $x, y, z$
with $x + y + z  \leq n$, for $n = 3, 4, \ldots $ .
This way, we recursively enumerate all solutions of
Fermat's equation. Hence, such a process
will eventually give a counterexample to Fermat's conjecture
\it if one exists%
\rm , but
the process will never yield conclusive evidence
if the conjecture happens to be true.
\end{comment}
In his famous address to the International
Mathematical Congress in 1900,
D. Hilbert\index{Hilbert, D.} proposed twenty-three mathematical problems
as a program to direct the mathematical efforts
in the twentieth century. The tenth problem\index{Hilbert's Tenth Problem} asks
for an algorithm which, given an arbitrary Diophantine
equation, produces either an integer solution for
this equation or indicates that no such solution exists.
After a great deal of preliminary work by
other mathematicians,
the Russian mathematician
Yuri V. Matijasevich finally
showed that no such algorithm exists.
But suppose we weaken the problem as follows.
First, effectively enumerate all Diophantine equations, and consider
the characteristic sequence $\Delta = \Delta_1 \Delta_2  \ldots $,
\index{characteristic sequence!of recursively enumerable set}
\index{set!recursively enumerable}
defined by $\Delta_i = 1$ if the $i$th Diophantine
equation is solvable, and 0 otherwise.
Obviously, $C( \Delta_{1:n} )  \leq n + O(1)$.

There is an algorithm to decide
the solvability of the first $n$
Diophantine equations
given about $\log n$ bits extra information.
Namely, given the number $m  \leq n$ of soluble
equations in the first $n$ equations, we can recursively
enumerate solutions to the first $n$ equations
in the obvious way
until we have found $m$ solvable equations.
The remaining equations are unsolvable. (This is a particular
case of a more general Lemma due to J.M. Barzdins known as
Barzdins' Lemma, \cite{LV90}.)

A.N. Kolmogorov observed that this shows that the solubility of
the enumerated Diophantine equations is
interdependent in some way since
$C ( \Delta_{1:n} | n)  \leq \log n +c$, for some fixed constant $c$.
The compressibility of the
characteristic sequence means in fact that the solvability
of Diophantine equations is highly interdependent---it
is impossible for an even moderately random sequence of them to be solvable
and the remainder unsolvable.

G.J. Chaitin proposed to replace
the question of mere solubility by the question
of whether there are finitely many or infinitely
many different solutions. That is,
no matter how many solutions we find
for a given equation, by itself this can give
no information on the question to be decided.
It turns out that
the set of indices of the Diophantine equations
with infinitely many different solutions
is not recursively enumerable.
In particular, in the characteristic sequence
\index{characteristic sequence!random}
each initial segment of length $n$ has Kolmogorov complexity
of about $n$.
 
\begin{claim}
\rm
There is an
(exponential) Diophantine equation
\[ A(n, x_1 , x_2 , \ldots ,x_m ) = 0\] which has
only finitely many solutions $x_1 , x_2 , \ldots ,x_m$
if the $n$th bit of $\Omega$ is zero and which has
infinitely many solutions $x_1 , x_2 , \ldots ,x_m$
if the $n$th bit of $\Omega$ is one.
\end{claim}
\begin{comment}
The role of {\em exponential} Diophantine equations
should be clarified.
Yu.V. Matijasevich\index{Matijasevich,
Yu.V.} [{\em Soviet Math.\ Dokl.},
11(1970), 354-357] proved that every recursively enumerable set has a polynomial
Diophantine representation.
J.P. Jones\index{Jones, J.P.}
and Yu.V. Matijasevich\index{Matijasevich, Yu.V.}
\index{Theorem!Jones-Matijasevich}
[{\it J.\ Symbol.\ Logic},
49(1984), 818-829]
proved that every recursively enumerable set 
has a singlefold exponential Diophantine
representation.
It is not known whether singlefold
representation (which is important in our
application) is always possible without exponentiation.
See also G.J. Chaitin\index{Chaitin, G.J.}, %
\it Algorithmic Information Theory%
\rm ,
Cambridge University Press, 1987.
\end{comment}
\begin{proof}
By dovetailing the running of
all programs of the reference prefix machine $U$ in the obvious
way, we find a recursive sequence of rational numbers
$\omega_1   \leq \omega_2    \leq \cdots$
such that $\Omega = \lim_{{n}   \rightarrow   \infty} \omega_n$.
The set
$$
R =   \{  (n, k): {\rm \ the\ }n {\rm th\ bit\ of\ }\omega_k
{\rm \ is\ a\ one}  \}
$$
is a recursively enumerable (even recursive) set.
The main step is to use
a theorem due to J.P. Jones\index{Jones, J.P.}
and Yu.V. Matijasevich\index{Matijasevich, Yu.V.}
\index{Theorem!Jones-Matijasevich}
[%
\it J.\ Symbol.\ Logic %
\rm 49(1984), 818-829] to the
effect that `every recursively enumerable set $R$ has a
singlefold exponential Diophantine representation $A( \cdot ,  \cdot )$'.
That is, $A(p, y) = 0$ is an exponential Diophantine equation,
and the singlefoldedness consists in the
property that $p   \in   R$
iff there is a $y$ such that $A(p, y) = 0$ is satisfied,
and, moreover, there is only a single such $y$.
(Here both $p$ and $y$ can be multituples
of integers, in our case $p$ represents $  \langle   n, x_1   \rangle   $, and
$y$ represents $  \langle    x_2 , \ldots ,x_m   \rangle   $.
For technical reasons we consider
as proper solutions only solutions $x$ involving no negative
integers.)
It follows that there is an exponential Diophantine equation
\sloppy{$A(n, k , x_2 , \ldots ,x_m ) = 0$ which has}
exactly one solution $x_2 , \ldots ,x_m$ if
the $n$th bit of the binary expansion of $\omega_k$
is a one, and it has no solution $x_2 , \ldots ,x_m$
otherwise.
Consequently, the number of different
$m$-tuples $x_1 , x_2 , \ldots ,x_m$ which are
solutions to $A(n, x_1 , x_2 , \ldots ,x_m ) = 0$
is infinite if the $n$th bit of the binary
expansion of $\Omega$ is a one, and this number is finite
otherwise.
\end{proof}

\end{example}

\subsection{Randomness of Individual Sequences Resolved}
The notion of randomness of an
infinite sequence in the sense of Martin-L\"of, as posessing
all effectively testable properties of randomness
(one of which is unpredictability), has turned out to be identical
with the notion of an infinite sequence having
the prefix Kolmogorov complexity of all finite initial
segments of at least the length of the initial segment itself,
Theorem~\ref{K3}.
This equivalence of a single notion being
defined by two completely different approaches is a truly remarkable fact.
(To be precise, the so-called  prefix Kolmogorov complexity
of each initial segment of the infinite binary sequence 
must not decrease more than a fixed constant, depending 
only on the infinite sequence,
below the length of that initial segment, \cite{LV90}.)
This property sharply distinguishes the random infinite binary
sequences from the nonrandom ones. The set of random infinite
binary sequences has uniform measure one. That means that as the outcome
from independent flips of a fair coin they occur with 
probability one.

For finite binary sequences the distinction between randomness
and nonrandomness cannot be abrupt, but must be a matter
of degree. For example, it would not be reasonable
if one string is random but becomes nonrandom if
we flip the first nonzero bit. In this context too
it has been shown that finite binary sequences which are
random in Martin-L\"of's sense correspond to those
sequences which have Kolmogorov complexity at least 
their own length. Space limitations forbid a complete
treatment of these matters here. Fortunately, it can be found
elsewhere, \cite{LV90}.

\section{Applications}
\subsection{Prediction}
We are given an initial segment of an infinite sequence
of zeros and ones. Our task is to predict the next
element in the sequence: zero or one? The set of possible
sequences we are dealing with constitutes the
`sample space'\index{sample space}; in this case, the set
of one-way infinite binary sequences. We assume some
probability distribution $\mu$ over the sample space,
where $\mu (x)$ is the probability of the initial segment
of a sequence being $x$. Then the probability of the next
bit being `0', after an initial segment $x$, is clearly
$\mu (0|x)=\mu (x0)/ \mu (x)$. This problem constitutes,
perhaps, the central task of
inductive inference and artificial intelligence.
However, the problem of induction
is that in general we do not know the distribution $\mu$,
preventing us from assessing the actual probability.
Hence, we have to use an estimate.

Now assume that $\mu$ is computable.
(This is not very restrictive, since any
distribution used in statistics is computable, provided the
parameters are computable.)
We can use
Kolmogorov complexity to give a very good
estimate of $\mu$. This involves the so-called
`universal distribution'\index{distribution!universal} ${\bf M}$.
Roughly speaking, ${\bf M} (x)$ is close to
$2^{-l}$, where $l$ is the length in bits of the shortest
effective description of $x$. Among other things, ${\bf M} $
has the property that
it assigns at least as high a probability to $x$ as
any computable $\mu$ (up to a multiplicative constant
factor depending on $\mu$ but not on $x$). What is
particularly important to prediction is the following.

Let $S_n$ denote the $\mu$-expectation of the square
of the error we make in estimating the probability
of the $n$-th symbol by ${\bf M}$. Then it can be shown
that the sum $\sum_n   S_n$ is bounded by a constant.
In other words, $S_n$ converges to zero
faster than $1/n$. Consequently, any actual (computable)
distribution can be estimated and predicted with
great accuracy using only the single universal distribution.
This approach is due to Raymond J Solomonoff (1926--- ), 
\cite{So60,So64,So78,So95} and predates 
Kolmogorov's complexity invention. This approach has led
to several developments in inductive reasoning,
including the widely applied celebrated 
`minimum message length' of Chris S. Wallace (1933--- ) 
and `minimum description length (MDL)' of Jorma Rissanen (1932--- ), 
[Chris S. Wallace and David M. Boulton, 
An information measure for classification,
{\em Computing Journal}, {\bf 11}(1968), 185-195; Jorma Rissanen,
Modeling by the shortest data description,
{\em Automatica-J.IFAC}, 14(1978), 465--471;
Jorma Rissanen, {\em Stochastical Complexity and Statistical Inquiry},
World Scientific Publishing Company, Singapore, 1989]
model selection methods
in statistics, by the analysis in
\cite{LV90}.

\subsection{G\"odel's incompleteness result}
\label{example.goedel}
\rm
\index{Theorem!Incompleteness}
We say that a formal system (definitions, axioms, rules of inference)
is {\em consistent}
if no statement which can be expressed in the system
can be proved to be both true and false in the system.
A formal system is {\em sound} if only true statements can be
proved to be true in the system. (Hence, a sound formal
system is consistent.)

Let $x$ be a finite binary
string. We write `$x$ is random' if
the shortest binary description of $x$ with respect to
the optimal specification method $D_0$ has length
at least $x$. A simple counting argument shows
that there are random $x$'s of each length.

Fix any sound
formal system $F$ in which we can express statements
like ``$x$ is random''. Suppose
$F$ can be described in $f$ bits---assume, for example,
that this is the
number of bits used in the exhaustive description of $F$
in the first chapter of the textbook %
\it Foundations of $F$%
\rm .
We claim that, for all but finitely
many random strings $x$, the sentence
`$x$ is random' is not provable in $F$.
Assume the contrary. Then given $F$, we can start to
exhaustively search for a proof that some string of length
$n \gg f $ is random, and print it when we find such a string $x$.
This procedure to print $x$ of length $n$ uses only $ \log  n  +  f$
bits of data, which is much less than $n$.
But $x$ is random by the proof and the fact that $F$ is sound.
Hence, $F$ is not consistent, which is a contradiction. This type of
argument is due to Janis M. Barzdins (1937--- ) 
and later G.J. Chaitin, see \cite{LV90}.

This shows that although most strings are random, it is impossible to
effectively prove them random. In a way, this explains why
the incompressibility method
is so successful. We can argue about a `typical' individual element,
which is difficult or impossible by other methods.

\subsection{Lower bounds}
\rm
The secret of the successful use of descriptional complexity
arguments as a proof technique is due to a simple
fact: the overwhelming majority of strings
have almost no computable regularities.
We have called such a string `random'.
There is no shorter description of such a string
than the literal description:
it is incompressible.
Incompressibility is a noneffective property in the sense of
Example~\ref{example.goedel}.

Traditional proofs often involve all instances
of a problem in order to conclude that some
property holds for at least one instance.
The proof would have proceeded simpler,
if only that one instance could have been
used in the first place.
Unfortunately, that instance is hard
or impossible to find, and the proof has to involve all the
instances.
In contrast, in a proof by the incompressibility method,
we first
choose a random (that is, incompressible) individual
object that is known to %
exist %
\rm (even though we cannot construct it).
Then we show that if the assumed property
would not hold, then this object could be compressed,
and hence it would not be random. Let us give a simple
example appearing in \cite{LV91,LV90}. 
A proof using  the probabilistic method
appears in 
[Paul Erd\"{o}s and Joel Spencer,
{\em  Probabilistic Methods in Combinatorics},
Academic Press, New York, 1974].

A {\it tournament}\index{tournament|bold}
is defined to be a complete directed graph. That is, for
\index{upper bounds!tournaments}
each pair of nodes $i$ and $j$, exactly one of edges
$(i,j)$ or $(j,i)$ is in $T$. The nodes of a tournament can be
viewed as players in a game tournament. If $(i,j)$ is in
$T$, we say player $j$ dominates player $i$.
We call $T$ {\it transitive}\index{tournament!transitive}
if $(i,j),(j,k)$ in $T$ implies $(i,k)$ in $T$.

Let $\Gamma = \Gamma_n$ be the set of
all tournaments on $N= \{ 1,  \ldots  , n \}$.
Given a tournament $T \in \Gamma$, fix a standard
encoding $E: T \rightarrow \{0,1\}^{n(n-1)/2}$,
one bit for each edge.
The bit for edge $(i,j)$ is set to 1 if $i < j$ and 0 otherwise.
There is a 1-1 correspondence between the members of $\Gamma$ and
the binary strings of length $n(n-1)/2$.

Let $v(n)$ be the largest integer such that every tournament
on $N$ contains a transitive subtournament on $v(n)$
nodes.

\begin{theorem}\label{tournament}
$v(n)  \leq  1 + \lfloor 2 \log n \rfloor $.
\end{theorem}
\begin{proof}
Fix $T \in \Gamma $ such that
\begin{equation}
C(E(T)|n,p) \geq n(n-1)/2,
\label{transrandom}
\end{equation}
where $p$ is a fixed program
that on input $n$ and $E'(T)$ (below)
outputs $E(T)$.
Let $S$ be the transitive subtournament of $T$
on $v(n)$ nodes. We try to compress $E(T)$,
to an encoding $E'(T)$, as follows.
\begin{enumerate}
\item
Prefix the list of nodes in $S$ in
order of dominance to $E(T)$,
each node using $\lceil \log n \rceil$ bits, adding
$v(n) \lceil  \log n \rceil$ bits.
\item
Delete all redundant bits from the $E(T)$ part, representing the
edges between nodes in $S$, saving $v(n)(v(n)-1)/2$ bits.
\end{enumerate}
Then,
\begin{equation}
l(E'(T)) = l(E(T)) -  \frac{v(n)}{2} (v(n)-1- 2 \lceil \log n \rceil ) .
\label{savings}
\end{equation}
Given $n$, the program $p$ reconstructs $E(T)$ from
$E'(T)$. Therefore,
\begin{equation}
C(E(T)|n,p) \leq l(E'(T)).
\label{compression}
\end{equation}
Equations~\ref{transrandom}, \ref{savings},
and \ref{compression} can only be satisfied with
$v(n) \leq 1+  \lfloor 2 \log n \rfloor$.
\end{proof}
The general idea used in the incompressibility
proof of Theorem~\ref{tournament} is the following. If each tournament
contains a large transitive subtournament, or any other
`regular' property for that matter, then also a tournament $T$ of
maximal complexity contains one. But the regularity induced by a too
large transitive subtournament can be used to compress the description
of $T$ to below its complexity, leading to the required
contradiction.

Results using the
incompressibility method can (perhaps) always
be rewritten using the
probabilistic method or counting arguments.
The incompressibility argument seems simpler and
more intuitive.
It is easy to generalize the above arguments from
proving `existence' to proving `almost all'.
Almost all strings have high complexity.
Therefore, almost all tournaments and
almost all undirected graphs have high complexity.
Any combinatorial property proven
about an arbitrary complex object in
such a class
will hold
for almost all objects in the class.
That is, such properties are subject to a Kolmogorov complexity 
0---1 Law: they either hold by a Kolmogorov complexity argument
for almost all objects in the class or for no objects in the class. 
For example, the proof of Theorem~\ref{tournament}
can trivially be strengthened as below. 
By simply counting the number of binary programs
of length at most the righthand side of Equation~\ref{random.exp}, 
there are at least
$2^{n(n-1)/2} (1 - 1/ n )$ tournaments $T$ on $n$ nodes
with
\begin{equation}\label{random.exp}
C(E(T)|n,p) \geq n(n-1)/2 -  \log n.
\end{equation}
This is a $(1-1/n)$th fraction
of all tournaments on $n$ nodes. Using Equation~\ref{random.exp}
in the proof
yields the statement below.

\begin{quote}
For almost all tournaments
on $n$ nodes (at least a $(1-1/n)$th fraction), the largest
transitive subtournament has at most $1+ 2  \lceil 2 \log n \rceil $
nodes, from some $n$ onwards.
\end{quote}

In \cite{LV90}, Chapter 6, the incompressibility
method is explained in more detail. Its utility is demonstrated in
a variety of examples of proving
mathematical and computational results.
These include questions concerning the
average case analysis of algorithms (such as Heapsort),
sequence analysis, average case complexity in general,
formal languages, combinatorics, time and
space complexity analysis
of various sequential or parallel machine models,
language recognition, and string matching. Other
applications include the use of resource-bounded Kolmogorov complexity
in the analysis of computational complexity classes,
the universal optimal search algorithm, and
`logical depth'.

\subsection{Statistical Properties of Finite Sequences}
Each individual infinite sequence generated
by a $( \frac{1}{2} , \frac{1}{2} )$
Bernoulli process (flipping a fair coin) has
(with probability 1)
the property that
the relative frequency of zeros in an initial $n$-length segment goes
to $\frac{1}{2}$ for $n$ goes to infinity.
Such randomness
related statistical properties of
individual (high) complexity finite binary sequences
are often required in applications of incompressibility
arguments.
The situation for infinite random sequences is studied above.

We know from Section~\ref{sect.random.infinite} that each
infinite sequence\index{sequence!random}
which is random with respect to the uniform measure
satisfies all effectively testable properties of randomness: it
is normal, it satisfies the so-called Law of the Iterated Logarithm,
the number of 1's minus the number of 0's in an initial
$n$-length segment is positive for infinitely many $n$ and
negative for another infinitely many $n$, and so on.
While the statistical properties of infinite sequences
are simple corollaries of the theory of Martin-L\"of randomness,
for finite sequences the situation is less simple.

In the finite case, randomness is a matter of degree,
because it would be clearly unreasonable to say
that a sequence $x$ of length $n$
is random and to say that a sequence $y$
obtained by flipping the first `1' bit of $x$ is nonrandom.
What we can do is to express the degree of
incompressibility of a finite sequence
in the form of its Kolmogorov complexity, and then
analyze the statistical properties of the sequence---for example,
the number of 0's and 1's in it, as in \cite{LV90,LV94}.

Since almost all finite sequences
have about maximal Kolmogorov complexity, each statistical property 
which is possessed
by each individual
maximal complexity sequence must also
hold approximately
in the expecteded sense (on average) for the overall set.
Let us look at the converse.
The fact that some property holds on the average
over the members of a set does not in general imply
that they are present in most individual members.
For example, if the set is $\{00\ldots 0, 11 \ldots 1\}$ where
both sequences have length $m$,
then the average 
relative frequency of 1's over initial segments of
length $n$ is $1/2$ for each $n$ with $1 \leq n \leq m$. 
But for each individual member
of this set the relative frequency of 1's in the
initial segments is either zero
or one.

In contrast, for infinite sequences all 
effectively testable properties
of randomness hold with probability one 
and therefore are also expected or average properties.
We have seen  that these properties
must also hold with certainty for each individual
sequence with high enough Kolmogorov complexity, and that
these sequences have probability one in the set
of all sequences. That is, the Kolmogorov random elements
in given set with uniform probability distribution
possess individually each effectively testable property which
is shared by almost all elements of the set.
We cannot expect that this phenomena also holds with
this sharpness in for the finite sequences,
if only because we cannot devide the set of finite sequences
sharply into sequences which are random and sequences
which are not random (unlike the infinite sequences).
However, it turns out that randomness properties hold
approximately for finite sequences with high enough 
Kolmogorov complexity.

For example, we can {\em a priori}
state that
each high-complexity finite binary sequence
is `normal' in the sense that each binary block of length
$k$ occurs about equally frequent for $k$ relatively small.
In particular, this holds for $k=1$.
However, in many applications we need to
know exactly what `about'
and the `relatively small' in this statement mean.
In other words, we are interested
in the extent to which `Borel normality' holds
in relation with the complexity of a finite sequence.

Let $x$ have length $n$.
It can be shown (and follows from the stronger result
below) that if $C(x|n) = n+O(1)$,
then the number of zeros it contains is
\[ \frac{n}{2} + O( \sqrt n ). \]

\begin{notation}
\rm
The quantity $K(x|y)$ in this section satisfies
\[ C(x|y) \leq K(x|y) \leq C(x|y)+2 \log C(x|y)+1.\]
This is
the length of
a self-delimiting version of a
program $p$ of length $l(p)=C(x|y)$, what we defined as
`prefix complexity'.
\end{notation}

\begin{definition}\label{def.normal}
\rm
The class of {\em deficiency} functions
is the set of functions
$\delta: {\cal N} \rightarrow {\cal N}$ satisfying
$K(n, \delta (n)|n-\delta(n))
\leq c_1$ for all $n$.
(Hence, $C (n, \delta (n)|n-\delta(n))
\leq c_1$ for all $n$.)
\end{definition}

This way we can retrieve $n$ and $\delta(n)$ from
$n- \delta (n)$ by a self-delimiting program
of at most $c_1$ bits.
We choose $c_1$ so large that each monotone
sublinear recursive function that we
are interested in, such as
$\log n$, $\sqrt n$, $\log \log n$, is such a deficiency function.
The constant $c_1$ is a benchmark which stays fixed throughout
this section.

\begin{lemma}
There is a constant $c$, such that for all deficiency functions
$\delta$, for each $n$ and $x \in \{0,1\}^n$, if
$ C(x) > n- \delta(n)$, then
\begin{equation}\label{result1}
\left| \#{\rm ones}(x)- \frac{n}{2} \right| < \sqrt{(\delta(n)+c)n \ln 2}.
\end{equation}
\end{lemma}

\begin{proof}
A general estimate of the tail probability of the binomial distribution,
with $s_n$ the number of successful outcomes in $n$ experiments
with probability of success $0 < p < 1$ and $q=1-p$,
is given by Chernoff's\index{Chernoff bounds}
bounds, \cite{LV90}
\begin{equation}\label{chernoffzerotex}
\Pr (|s_n -np| \geq m ) \leq 2e^{-m^2 /4npq } .
\end{equation}

Let $s_n$ be the number of 1's
in the outcome of $n$ fair coin flips, which
means that $p=q=1/2$.
Defining $A=\{ x \in \{0,1\}^n : |\#{\rm ones}(x) - n/2|\geq m \}$
and applying Equation~\ref{chernoffzerotex},
\[ d(A) \leq 2^{n+1}e^{-m^2 /n}. \]
Let $m= \sqrt{(\delta(n) + c )n\ln 2}$ where $c$ is a constant
to be determined later.
We can compress any $x \in A$ in the following way.
\begin{enumerate}
\item
Let $s$ be a self-delimiting program to retrieve $n$ and $\delta (n)$
from $n-\delta(n)$, of
length at most $c_1$.
\item
Given $n$ and $\delta(n)$, we can effectively enumerate $A$.
Let $i$ be the index of $x$ in such an
effective enumeration of $A$.
The length of the (not necessarily self-delimiting) description of $i$
satisfies

\begin{eqnarray*}
l(i) \leq \log d(A) & = & n+1 + \log e^{-m^2 /n}\\
 & \leq & n + 1 - \delta(n) - c.
\end{eqnarray*}
\end{enumerate}

The string $si$ is padded to length $n+1 - \delta(n) - c + c_1$.
From $si$ we can reconstruct $x$ by first using $l(si)$ to compute $n
- \delta(n)$, then compute $n$ and $\delta(n)$
from $s$ and $n-\delta(n)$, and subsequently
enumerate $A$ to obtain the $i$th element. Let $T$ be the Turing
machine embodying the
procedure for reconstructing $x$. Then, by definition of $C(\cdot)$,
\[ C(x) \leq C_{T}(x) + c_{T} \leq n+1 - \delta(n)
- c + c_1 + c_{T}. \]
Choosing $c = 1+c_1 + c_{T}$ we find
$C(x) \leq n- \delta (n)$,
which contradicts the condition of the theorem.
Hence, $ |\#{\rm ones}(x)- n/2 | <m$.
\end{proof}

It may be surprising at first glance, but there
are no maximally complex sequences with about equal
number of zeros and ones. Equal numbers of zeros and ones is
a form of regularity, and therefore lack of complexity.
That is, for $x \in \{0,1\}^n$, if $|\#{\rm ones}(x) - n/2| = O(1)$,
then the randomness deficiency $\delta(n)=n-C(x)$ is
nonconstant (order $\log n$).

The analysis up till now has been about the statistics
of 0's and 1's. But in a normal
infinite binary sequence according to Definition~\ref{def.normal}
each block of length $k$ occurs with limiting frequency of $2^{-k}$.
That is, blocks 00, 01, 10, and 11 should
occur about equally often, and so on. Finite sequences
will generally not be exactly normal, but normality
will be a matter of degree.
We investigate the block statistics for finite binary
sequences.

\begin{definition}\label{def.xynl}
\rm
Let $x=x_1 \ldots x_n$ be a binary string of
length $n$, and $y$ a much smaller string
of length $l$. Let $p = 2^{-l}$ and
$\#y(x)$ be the number of
(possibly overlapping) distinct occurrences of $y$ in $x$.
For convenience, we assume that
$x$ `wraps around' so that an occurrence
of $y$ starting at the end of $x$
and continuing at the start also counts.
\end{definition}

\begin{theorem}\label{blockszerotex}
Assume the notation of Definition~\ref{def.xynl} with $l \leq \log n$.
There is a constant $c$,
such that for all $n$ and $x \in \{0,1\}^n$,
if $C(x) > n - \delta(n)$, then
\[ |\#y(x)-np|  <  \sqrt{ \alpha np},\]
with $\alpha = [K(y|n)+\log l + \delta(n)+c](1-p)l 4 \ln 2 $.
\end{theorem}

It is known from probability theory
that in a randomly generated finite sequence the
{\em expectation} of the length of the longest run of zeros or ones
is pretty high. For each individual finite sequence with high Kolmogorov
complexity we are {\em certain} that it contains each
block (say, a run of zeros) up to a certain length.
\begin{theorem}\label{zeros}
Let $x$ of length $n$ satisfy $C(x) \geq n-\delta(n)$.
Then, for sufficiently large $n$, $x$ contains all blocks $y$ of length

\[ l= \log n - \log \log n - \log (\delta(n) + \log n) - O(1). \]
\end{theorem}
\begin{corollary}
\rm
If $\delta(n)=O(\log n)$, then each
block of length $\log n - 2 \log \log n -O(1)$
is contained in $x$.
\end{corollary}
Analyzing the proof of Theorem~\ref{zeros} we can
improve this when $K(y|n)$ is low.
\begin{corollary}
\rm
If $\delta(n)=O( \log \log n)$, then for each $\epsilon > 0$
and $n$ large enough, $x$ contains an all-zero run $y$
(for which $K(y|n)=O(\log l)$)
of length $l= \log n - (1+ \epsilon ) \log \log n +O(1)$.
Since at least a fraction of $1- 1/\log n$ of all strings of length
$n$ has such a $\delta (n)$, the result gives the {\em expected}
value of the length of the longest all-zero run.
This is almost a $\log \log n$ additional term
better than the expectation reported using simple
probabilistic methods in 
[Thomas H. Cormen, Charles E. Leiserson, Ronald L. Rivest,
{\em Introduction to Algorithms}, MIT Press, Cambridge, Mass., 1990].
\label{expect}
\end{corollary}

\subsection{Chaos and Predictability}
\rm
Given sufficient information about a physical system,
like the positions, masses and velocities of all particles,
and a sufficiently powerful computer
with enough memory and computation time,
it should be possible in principle to compute
all of the past and all of the future of the system.
This view, eloquently propagated by
P.S. Laplace\index{Laplace, P.S.}, can be espoused
both in classical mechanics and quantum mechanics.
In classical mechanics one would talk about a
single `history', while in quantum mechanics
one would talk about probability distributions over
an ensemble of `possible histories.'
Nonetheless, in practice it is impossible to
obtain all parameters precisely. The finitary nature
of measurement and computation requires truncation of
real valued parameters; there are measuring errors;
and according to basic quantum mechanics it is
impossible to measure certain types of parameters
simultaneously and precisely. Altogether, it is fundamental
that there are minute uncertainties in our knowledge
of any physical system at any time.

This effect can been combined with the consistent tradition that
small causes can have large effects exemplified by
the metaphor
``a butterfly moving its wing in tropical
Africa can eventually cause a cyclone in the Caribbean.''
Minute perturbations in initial conditions can
cause, mediated by strictly computable
functions, arbitrary large deviations in outcome.
In the mathematics of nonlinear deterministic systems
this phenomenon has been described by the catch term `chaos'.

The unpredictability of this phenomenon is
sometimes explained through Kolmogorov complexity.
Let us look at an example where
unpredictability is immediate (without using
Kolmogorov complexity).

The following mathematical conditions
are satisfied by classic thermodynamic systems.
Each point in a state space
$X$ describes a possible
micro state of the system. Starting from each position,
the point representing the system describes a trajectory
in $X$. We take time to be discrete and
let $t$ range over the integers.

If the evolution of the system in time
can be described by a transformation group $U^t$
of $X$
such that $U^t  \omega $ is the trajectory starting
from some point $\omega$ at time 0, then we call
the system {\em dynamical}.
Assume that $U^t \omega$ is computable as a function
of $t$ and $\omega$ for both positive and
negative $t$.

\begin{definition}\label{def.physical.space}
\rm
Assume that the involved measure $\mu$ is recursive
and satisfies
$\mu (U^t V)= \mu (V)$
for all $t$ and all measurable sets $V$ (like the $\Gamma_x$'s),
as in Liouville's Theorem.
The $\mu$-measure
of a volume $V$ of points in state space is invariant over time.
A {\em physical system} ${\cal X}$ consists of the
space $X$, $n$-cells $\Gamma_x$,
transformation group $U^t$, and the recursive
and volume invariant measure $\mu$.
\end{definition}

Consider a system ${\cal X}$ 
with 
the {\em orbit} of the system 
the sequence of subsequent micro states
$\omega , U \omega , U^2 \omega , \ldots$.

We cannot measure states which are (or involve) real
numbers precisely. Because of the finite precision
of our measuring instruments, the
number of possible {\em observed states} is finite.
Hence, we can as well assume that the {\em observable state space}
${\cal S}$ of our system contains only a finite number of states,
and constitutes a finite partition of ${\cal X}$.
Then, each micro state in ${\cal X}$ is observed as
being an element of some observable state in ${\cal S}$.

For convenience we assume that
${\cal S}=\{ \Gamma_0, \Gamma_1\}$ is
the set of observed states. Here $\Gamma_0$
is the set of micro states starting with a `0',
and $\Gamma_1$ is the set of micro states starting with a `1'.

Let $\Gamma^n$ be the observed state of $X$ at time $n$.
Given an initial observed segment
$\Gamma^0 , \ldots , \Gamma^{n-1}$ of an observed orbit we want
to compute the next observable state
$\Gamma^{n}$. Even if it would not be possible
to compute it, we would like to compute a
prediction of it which does better than a random coin flip.

A well-known example of a chaotic system is the
the {\em doubling map} which results
from the so-called
baker's map
by deleting all bits with non-positive indexes from
each state. That is, each micro state $\omega$ is a
one-way infinite binary sequence. Considering
this as the binary expansion
of a real number after the binary point, the set
of micro states is the real interval $[0,1)$.
The system
evolves according to the transformation
\begin{equation}\label{ford.1}
\omega_{n+1} = 2\omega_n \pmod{1}
\end{equation}
where `mod 1' means drop the integer part. All iterates of
Equation~\ref{ford.1} lie in the unit interval $[0,1)$.
The observable states are
$\Gamma_0 = [0, \frac{1}{2} )$ and
$\Gamma_1 =[\frac{1}{2} , 1)$.
\noindent
\begin{comment}
\rm
The doubling map is related to
the discrete logistic equation
\[ Y_{n+1} = \alpha Y_n (1 - Y_n ) \]
which maps the unit interval upon itself when $0 \leq \alpha \leq 4$.
When $\alpha =4$, setting $Y_n = \sin 2 \pi X_n$,
we get precisely
\[ X_{n+1} = 2 X_n \pmod{1} . \]
\end{comment}
Assuming that the initial state is randomly drawn from $[0,1)$
according to the uniform measure $\lambda$, we can use
complexity arguments to show that the doubling map's observable
orbit cannot be predicted better than a coin toss
[Joseph Ford, `How random is a random coin toss?',
{\em Physics Today}, {\bf 36}(1983),
April issue, pp. 40-47].

Namely, with $\lambda$-probability 1 the drawn
initial state will be a Martin-L\"of random
infinite sequence since we have shown that they have uniform measure
one in the set of infinite sequences. Such random sequences
by definition cannot be effectively predicted
better than
a random coin toss, \cite{LV90}.

But in this case we do not need to go
to such trouble. The observed orbit essentially
consists of the consecutive bits of the initial state.
Selecting that initial state randomly from the
uniform measure is isomorphic to flipping a fair
coin to generate it.

\end{document}